\documentclass{amsart}
\usepackage{graphicx}
\usepackage[utf8]{inputenc}
\usepackage[T1]{fontenc}
\usepackage{amsmath}
\usepackage{amsfonts}
\usepackage{amssymb}
\usepackage[version=4]{mhchem}
\usepackage{stmaryrd}
\usepackage{mathrsfs}
\usepackage{amssymb,amscd,amsthm,amsxtra}
\usepackage{latexsym}
\usepackage{epsfig}
\usepackage{color}
\usepackage{cases}
\usepackage{esint}
\usepackage{mathtools}
\usepackage{empheq}
\usepackage{pgfplots}
\pgfplotsset{compat=1.17}
\usepackage{tikz}
\usepackage{xcolor}
\usepackage[most]{tcolorbox}
\usetikzlibrary{arrows.meta}
\newtheorem{thm}{Theorem}[section]
\newtheorem{cor}[thm]{Corollary}
\newtheorem{lem}[thm]{Lemma}
\newtheorem{prop}[thm]{Proposition}
\theoremstyle{definition}
\newtheorem{defn}[thm]{Definition}
\theoremstyle{remark}
\newtheorem{rem}[thm]{\bf Remark}
\numberwithin{equation}{section}

\DeclareMathOperator*{\esssup}{ess \sup}

\newcommand{\set}[1]{\left\{#1\right\}}
\newcommand{\R}{\mathbb R}
\newcommand{\be}{\begin{equation}}
\newcommand{\ee}{\end{equation}}

\newcommand{\eps}{\varepsilon}

\newcommand{\Pa}{\mathcal{P}}
\newcommand{\Qt}{\tilde Q}

\begin{document}

\title[Phragm\'en-Lindel\"of theorems for equations of $p$-laplacian-type]{Phragm\'en-Lindel\"of theorems for parabolic equations of $p$-laplacian-type}
\author[Enrico Zani]{Enrico Zani} 
\address{Dipartimento di Matematica ``F. Casorati'', Universit\`a di Pavia, Via Ferrata, 5, 27100 Pavia, Italy.}
\email{\tt e.zani12@campus.unimib.it}


\begin{abstract}
We consider rigidity properties of ancient solutions to parabolic equations of $p$-Laplacian-type. We start with the degenerate range $p>2$. First, a solution which is bounded above in an intrinsic backward paraboloid is either bounded in an entire half-space, or grows in space at least like a power with a sharp, explicit exponent. This relies on the expansion of positivity and on the Comparison Principle. Then, we turn to the growth in time. Positivity in the paraboloid alone yields no upper rate: we present explicit self-similar solutions, positive in the paraboloid, with $u(0,-T)\simeq T^{\beta}$ for every $\beta\le0$ and every $\beta\ge1$. For solutions bounded below in the paraboloid, we prove a dichotomy for the infimum of $u(\cdot,-T)$ on a ball whose radius is a fixed fraction of that of the paraboloid at time $-T$: either it is of order $T$ along a sequence of times, or it grows sublinearly, like $T^{\zeta}$ with $\zeta<1$ explicit. In the second case, a Harnack-type argument shows that the supremum on the same ball grows at most linearly. Then, we carry out the same program in the singular, supercritical range $\frac{2N}{N+1}<p<2$, where the picture is similar but several roles are reversed. In space, the iteration via the Comparison Principle can only be started from a point lying above an explicit slope, and explicit solutions show that nothing better than a cone bound can be expected below it. The self-similar profiles are unique but only local, and positivity in the paraboloid is compatible with every exponent $\beta\le1$ and with no exponent in $(1,\frac1{2-p})$. In time, the dichotomy is replaced by an a priori linear bound for the infimum. Similarities and differences between the two ranges are collected at the end of the work.
\end{abstract}

\maketitle
\begin{center}
		\begin{minipage}{9cm}
			\small
			\tableofcontents
		\end{minipage}
	\end{center}

\section{Introduction}
The aim of this work is to present a new result of Phragm\'en-Lindel\"of-type concerning ancient solutions to parabolic $p$-Laplace equations, such as 
\begin{equation}
\label{toy}
    u_t - \operatorname{div}(|Du|^{p-2} Du) =0.
\end{equation}
The literature about parabolic $p$-Laplace equations counts a huge number of works. A decisive contribution is due to DiBenedetto, with a series of works such as \cite{DPE} and \cite{DGV}, just to cite two of the most important ones. There, the issues about the existence and the regularity of the solutions are treated in detail. 

Concerning such equations, rigidity results, such as Liouville theorems, do exist. Indeed, such results can be seen as consequences of the Harnack inequality, which holds both for linear and nonlinear equations. For instance, the linear case with continuous coefficients can be seen in \cite{Landis}, and \cite{Glagoleva1} also considered discontinuous coefficients. 
Similarly to the nondegenerate case, for which bounded solutions in $t<0$ are constant, \cite{DGV-Liouville} and \cite{DGV} prove that every (two-sided) bounded solution to \eqref{toy} in $\R^{N} \times (-\infty,T)$ is constant. 

It is then worth asking whether contiguous rigidity results, such as theorems of Phragm\'en-Lindel\"of-type, do hold. These families of theorems are different in the conclusions, but are similar in the assumptions. Indeed, both of them start by proper hypotheses on the sign of the solution. While Liouville-type results aim to prove the constancy of $u$, Phragm\'en-Lindel\"of theorems try to infer an optimal growth rate, with respect to space or time.

The first Phragm\'en-Lindel\"of theorem, \cite{PL}, characterized holomorphic and harmonic functions in an unbounded planar domain, resulting in a superlinear growth with respect to the space. Afterwards, a consistent literature proved theorems of this kind not only for more general elliptic equations, but also for parabolic equations, dealing with the growth in time as well. For instance, the linear case was treated in works as \cite{Shi} and \cite{Gla}. While the linear equation seems to be settled, this is not the case for the degenerate nonlinear equation \eqref{toy}, at least in $\R^{N} \times (-\infty,T)$. 

\subsection{Novelty and significance}
While for linear equations the picture seems clear (see \cite{Shi} and \cite{Gla}), to the best of our knowledge no analogue was available for the degenerate equation. For $p>2$ the paper contributes in three directions, one in space and two in time; the same three issues are then revisited in the singular, supercritical range $\frac{2N}{N+1}<p<2$. In this work, \emph{supercritical} is meant as in \cite{DGV}.

\smallskip\noindent\emph{Optimal growth in space.} Theorem~\ref{P-L Thm in space} gives a growth rate in space for solutions which are bounded above in the paraboloid $\Pa_h$. This is a remarkable geometry: the boundary of $\Pa_h$ is the zero set of an explicit eternal solution (Remark~\ref{funzione paraboloidica}), so the paraboloid is singled out by the equation itself, and it is the natural counterpart of the paraboloid used by Glagoleva in the linear case, \cite{Gla}. In this setting we prove an optimal result: either the solution is bounded in a half-space, or it grows at least like $|x|^{\frac p{p-1}}$, and an explicit solution shows that this exponent is attained. The alternative therefore characterizes completely the problem in space for this geometry. The proof only uses the expansion of positivity and the Comparison Principle, so it holds for the whole class of operators satisfying the structure conditions of Section~\ref{sec:equations}, and not only for the prototype.

\smallskip\noindent\emph{Explicit solutions and growth in time.} In time the situation is different: positivity in $\Pa_h$ alone yields no upper rate. In Section~\ref{section4} we consider explicit self-similar solutions, whose profile solves a nonlinear ordinary differential equation, and we determine for which exponents they are positive in $\Pa_h$. Self-similar solutions are a classical tool for the parabolic $p$-Laplacian and for the porous medium equation, see for instance \cite{KV} for the fundamental solutions and the asymptotic behavior of the former, and \cite{V} for a comprehensive treatment of the latter; here we use the ansatz backward in time, which leads to a different profile equation. Power growth in time of arbitrarily high order is compatible with positivity, whereas a whole range of exponents is excluded. This is the non-trivial outcome of the explicit solutions: a result in time cannot be an a priori bound, and it can only be a dichotomy.

\smallskip\noindent\emph{A dichotomy in time.} Section~\ref{section5} shows what remains true. By a DeGiorgi iteration, which uses only the structure conditions, the infimum on a fixed fraction of the slices is either comparable to the time along a sequence, or it grows at most like a power with an explicit exponent smaller than one; in the latter case, the Harnack inequality gives a linear bound for the supremum. The explicit solutions of Section~\ref{section4} realize both alternatives, so neither of them can be excluded.

\smallskip\noindent\emph{The singular supercritical range.} The same questions make sense for $\frac{2N}{N+1}<p<2$, where solutions are locally H\"older continuous and the Harnack inequality holds. In Sections~\ref{section6}--\ref{section9} we show that the same program survives, but with quite different outcomes. In space (Theorem~\ref{s:thm:space}) the iteration can only be started from a point lying above an explicit slope $\Xi_*$, because the time scale of the expansion of positivity now \emph{increases} with the size of the solution; the explicit solutions of Section~\ref{section7} show that below this slope only a cone bound is available. The self-similar profiles are unique, but they are only local when $\beta<0$, and the range of exponents compatible with positivity in $\Pa_h$ is $\beta\le1$, instead of $\beta\le0$ or $\beta\ge1$. In time (Theorem~\ref{s:thm:time}) the persistence lemma works in the opposite direction and the dichotomy becomes an a priori linear bound for the infimum. The step from the infimum to the supremum is not transposed. Section~\ref{section10} collects similarities and differences between the two ranges.

\subsection{Plan of the paper}
In Section~\ref{section2} we fix the notation and recall the structure conditions, the notion of weak solution, and the three tools we rely on: the expansion of positivity, the Harnack inequality (both stated for $p>2$ and for the singular supercritical range, since the statements differ) and the Comparison Principle. Sections~\ref{section3}--\ref{section5} deal with the degenerate case $p>2$. Section~\ref{section3} contains the result in space (Theorem~\ref{P-L Thm in space}) and the explicit solution showing that it is sharp. Section~\ref{section4} is devoted to the self-similar solutions and to the positivity criterion (Theorems~\ref{thm:profile} and~\ref{thm:positivity}). Section~\ref{section5} contains the dichotomy in time (Theorem~\ref{thm:dich}), the linear bound on the supremum (Theorem~\ref{thm:sup}), and a comparison with the explicit solutions of Section~\ref{section4}. Sections~\ref{section6}--\ref{section9} are the singular counterparts, for $\frac{2N}{N+1}<p<2$, in the same order: the result in space (Theorem~\ref{s:thm:space}) and the explicit solutions showing why the iteration cannot be started from every point (Section~\ref{section7}), the self-similar profiles and the positivity criterion (Proposition~\ref{s:prop:shape} and Theorem~\ref{s:thm:positivity}), and the a priori linear bound in time (Theorem~\ref{s:thm:time}), which replaces the dichotomy. Section~\ref{section10} summarizes similarities and differences between the two ranges (for instance, the same threshold $q_\mu$ but reversed inequalities, and the same positivity criterion but exchanged ranges of $\beta$).

\section{Preliminaries}
\label{section2}
In the following lines, we fix the notation and we recall some fundamental definitions and results, to be used in the following sections. For a wider and complete treatise and for the relative proofs, we refer to \cite{DGV}. 

\subsection{Notation}\label{sec:notation}
The notation has been chosen so as to dispel as many doubts as possible. In particular, subscripts are used to mark the differences between the choices of constants: a letter carrying a subscript has a fixed meaning, which is recalled where the letter is introduced. We collect here the elementary notation.
\begin{itemize}
\item The dimension is $N\ge1$. We write $|x|$ for the Euclidean norm of $x\in\R^N$ and $|x|_\infty:=\max_i|x_i|$. The Lebesgue measure of a measurable set $\Omega$ is denoted by $|\Omega|$, and $D$ is the gradient in the space variables. The gradient variable of the vector field $\mathbf A$ is always denoted by $\xi\in\R^N$.
\item For $y\in\R^N$ and $\rho>0$ we consider both the open ball $B_\rho(y):=\{x\in\R^N:\ |x-y|<\rho\}$ and the open cube of edge $2\rho$ centered at $y$, with faces parallel to the coordinate planes, $K_\rho(y):=\{x\in\R^N:\ |x-y|_\infty<\rho\}$. We write $B_\rho:=B_\rho(0)$ and $K_\rho:=K_\rho(0)$. Both centers and radii vary, and $B_\rho(y)\subset K_\rho(y)\subset B_{\sqrt N\rho}(y)$. The letter $K$ is used only for cubes.
\item $C$ denotes a generic positive constant. Its exact value is irrelevant, unless otherwise specified by a subscript or by the context, and it may change from line to line, even within the same chain of inequalities. When needed, we specify the dependencies: for instance, $C_{T_\star}$ in the proof of Theorem~\ref{thm:sup} is a constant depending on the time $T_\star$ (and on $u$ and $\mu$) fixed there, and it is not related to the structure constants $C_0,C_1$ of Section~\ref{sec:equations}, nor to the constant $c_1$ of \eqref{eq:mu0c1}. The expression ``depending only on the data'' means depending only on $p,N,C_0,C_1$. Constants with a fixed meaning always carry a subscript or a different letter, such as $C_0,C_1,C_H,\Gamma_H,\Gamma_1,\Gamma_2,\Gamma_3,c_0,c_1,\kappa_0$; the letter $c$ without subscript is used only for the amplitude of the self-similar solutions $u_{\beta,c,s}$ of Section~\ref{section4}. Similarly, $M$ denotes a positive number which may be large, and $M_\beta$, $M_\mu(T)$ have the specific meaning given where they are introduced.
\item Symbols introduced inside a proof (for instance $\xi_\eta,\omega,\Lambda_\eta$ in the proof of Theorem~\ref{P-L Thm in space}, or $E(r)$ in the proof of Proposition~\ref{prop:uniq}) have a local meaning, which ends with the proof.
\item In $E_T=E\times(0,T]$ (Section~\ref{sec:equations}) the number $T>0$ only gives the length of the cylinder; in Sections~\ref{section4},~\ref{section5},~\ref{section8} and~\ref{section9}, instead, $T>0$ is a positive time, and the corresponding time level is $-T$.
\item Roman numerals $\text{(i)},\text{(ii)}$ label the two alternatives of Theorem~\ref{P-L Thm in space}, of Theorem~\ref{s:thm:space} and of Theorem~\ref{thm:dich}. They are always accompanied by the reference to the relevant theorem, and conditions or cases appearing inside a proof are labelled differently.
\end{itemize}

\subsection{Equations of $p$-Laplacian-type}\label{sec:equations}
Let $E$ be an open set in $\R^{N}$, and for $T>0$ let $E_{T}$ stand for the cylindrical domain $E \times (0,T]$. 
Moreover, given a cylinder $Q\subset E_T$, by $\partial_p Q$ we mean its parabolic boundary. 
Given $p>1$, we consider homogeneous, quasilinear, parabolic partial differential equations of $p$-Laplacian-type as follows:  
\begin{equation}\label{p-Laplace}
    u_{t}-\operatorname{div} \mathbf{A}(x, t, Du)= 0 \text{ weakly in } E_T,
\end{equation}
where $\mathbf{A}:E_T \times \R^{N} \to \R^N$ is measurable, independent of $u$, odd in the gradient variable $\xi\in\R^N$, and satisfies the structure conditions 
\begin{equation}
\label{structure}
\begin{cases}
    \mathbf{A}(x, t, \xi) \cdot \xi \geq C_{0}|\xi|^{p}, \\
    |\mathbf{A}(x, t, \xi)| \leq C_{1}|\xi|^{p-1} 
\end{cases}    
\end{equation}
for a.e. $(x,t)\in E_T$ and every $\xi\in\R^N$, where $0<C_0\le C_1<\infty$. We also assume $\mathbf{A}$ to be monotone in the gradient variable, in the sense that 
\begin{equation}
\label{monotonicity}
(\mathbf{A}(x,t,\xi_1)-\mathbf{A}(x,t,\xi_2)) \cdot (\xi_1 - \xi_2) \geq 0,
\end{equation}
for every $(x,t)\in E_T$ and $\xi_1$, $\xi_2 \in \R^N$. The prototype of such a class of parabolic equations is 
\begin{equation}\label{Eq:prot}
    u_{t}-\operatorname{div}\left(|D u|^{p-2} Du\right)=0, \ \text{ weakly in } E_T.
\end{equation}
The equation \eqref{p-Laplace}-\eqref{monotonicity} is said to be {\it degenerate} if $p>2$, {\it singular} if $1<p<2$, and {\it nondegenerate} if $p=2$. 

Another example is the equation
\begin{equation}\label{Eq:prot2}
u_{t}-\big(a_{ij}(x, t)\,|D u|^{p-2}\, u_{x_i}\big)_{x_j}=0, \ \text { weakly in } E_{T},
\end{equation}
with the summation convention, where $(a_{ij}(x,t))$ is a measurable symmetric matrix
which is uniformly positive definite in $E_T$, that is, for some $C\ge1$,
\[
C^{-1}|\xi|^2\le a_{ij}(x,t)\,\xi_i\xi_j\le C|\xi|^2\qquad\text{for a.e. }(x,t)\in E_T\text{ and every }\xi\in\R^N.
\]
Here $\mathbf A(x,t,\xi)=|\xi|^{p-2}\,(a_{ij}(x,t)\xi_i)_{j=1,\dots,N}$ satisfies \eqref{structure}
with $C_0=C^{-1}$ and $C_1=C$. If $1<p<2$, $\mathbf A$ is extended by $\mathbf A(x,t,0)=0$; this is the continuous
extension, since $|\mathbf A(x,t,\xi)|\le C_1|\xi|^{p-1}\to0$ as $\xi\to 0$. As for \eqref{monotonicity}, a direct computation of the
Jacobian matrix of $\xi\mapsto|\xi|^{p-2}(a_{ij}\xi_i)_j$ shows that it holds, for instance,
if $a_{ij}=a\,\delta_{ij}$ for a scalar function $a$, or more generally if
$|p-2|\,C^2<p$, which reads $(p-2)\,C^2<p$ for $p>2$ and $(2-p)\,C^2<p$ for $1<p<2$
(and is void for $p=2$).

A further yet similar example is 
\[
\mathbf A(x,t,\xi)=\big(a_{ij}(x,t)\xi_i\xi_j\big)^{\frac{p-2}{2}}\,(a_{ij}(x,t)\xi_i)_{j=1,\dots,N}
= D_\xi\Big[\tfrac1p\big(a_{ij}(x,t)\xi_i\xi_j\big)^{\frac p2}\Big],
\]
extended by $\mathbf A(x,t,0)=0$ (which is the continuous extension, also when $1<p<2$).
It is odd in $\xi$, and it is the gradient of the convex $C^1$ function
$\xi\mapsto\frac1p\|\xi\|_{a}^p$, where $\|\xi\|_a=(a_{ij}\xi_i\xi_j)^{1/2}$ is a norm;
hence $\mathbf A$ is monotone as in \eqref{monotonicity}. Moreover \eqref{structure} holds
for every $p>1$ with $C_0=C^{-\frac p2}$ and
$C_1=C^{\frac p2}$ if $p\ge2$, $C_1=C^{\frac{4-p}{2}}$ if $1<p<2$.

\begin{defn}\label{def:weak}
A function $u \in C^{0}_{loc}\left(0, T; L_{loc}^{2}(E)\right) \cap L_{loc}^{p}(0, T ; W_{loc}^{1,p}(E))$ is a weak sub(super)-solution to \eqref{p-Laplace}-\eqref{monotonicity} if for every $[t_1,t_2] \subset (0,T]$ and for every nonnegative test function $\varphi \in W^{1,2}_{loc}(0,T; L^{2}(E)) \cap L_{loc}^{p}(0,T;W^{1,p}(E))$ such that $\operatorname{supp}(\varphi)\cap(\R^N\times[t_1,t_2])$ is a compact subset of $E\times[t_1,t_2]$,
    \[
    \left.\int_{E} u \varphi\, dx \right|_{t_1}^{t_2} + \int_{t_1}^{t_2} \int_{E} [-u \varphi_t + \mathbf{A}(x, t, Du) \cdot D\varphi] \ dx \ dt 
    \leq (\geq) 0 .\]
\end{defn}
In Sections~\ref{section3}--\ref{section5} we deal with degenerate equations, namely $p>2$, and in Sections~\ref{section6}--\ref{section9} with singular equations in the supercritical range $\frac{2N}{N+1}<p<2$. For both ranges \cite{DGV} proves in detail notable properties for solutions, such as the expansion of positivity, the Harnack inequality, and the comparison principle. We hereby recall them in their essential features; the first two have different statements in the two ranges, and are therefore stated separately.

The so-called expansion of positivity turns a measure-theoretical information in a pointwise estimate on the solution. 
\begin{prop}[Degenerate Expansion of Positivity]
\label{exp pos}
Let $u$ be a nonnegative, weak super-solution to \eqref{p-Laplace}-\eqref{monotonicity}. Given $M>0$, suppose that for some $(y,s) \in E_T$ and some $\rho>0$ 
\begin{equation}
\label{exp pos datum}
|[u(\cdot, s) \geq M] \cap K_{\rho}(y)| \geq \alpha |K_{\rho}(y)|    
\end{equation}
for some $\alpha \in (0,1)$. Then, there exist constants $\eta,$ $\delta \in (0,1)$ and $b>1$, depending only on the data and $\alpha$, such that 
\begin{equation}
u(\cdot, t) \geq \eta M \text{ a.e. in } K_{2\rho}(y)    
\end{equation}
for all times 
\begin{equation}
\label{exp pos time}
s+ \frac{1}{2} \frac{b^{p-2}}{(\eta M)^{p-2}}\delta \rho^p \leq t \leq s+ \frac{b^{p-2}}{(\eta M)^{p-2}}\delta \rho^p,    
\end{equation}
provided $K_{8\rho}(y) \times (s, s+ \frac{b^{p-2}}{(\eta M)^{p-2}}\delta \rho^p] \subset E_T$.
\end{prop}
We stress that $\eta,\delta,b$ are constants depending only on the data and $\alpha$: the symbol $\delta_\mu$ of Section~\ref{section5} denotes, instead, a number depending only on $\mu$ and $p$.

In the singular range the time scale of the expansion of positivity changes, since $M^{2-p}$ now \emph{increases} with $M$, and the cube $K_{8\rho}$ is replaced by $K_{16\rho}$.
\begin{prop}[Singular Expansion of Positivity]
\label{exp pos sing}
Let $1<p<2$ and let $u$ be a nonnegative, weak super-solution to \eqref{p-Laplace}-\eqref{monotonicity}. Given $M>0$, suppose that for some $(y,s) \in E_T$ and some $\rho>0$ the datum \eqref{exp pos datum} holds for some $\alpha \in (0,1)$. Then, there exist constants $\eta,$ $\delta,$ $\bar\varepsilon \in (0,1)$, depending only on the data and $\alpha$, such that
\begin{equation}
u(\cdot, t) \geq \eta M \text{ a.e. in } K_{2\rho}(y)
\end{equation}
for all times
\begin{equation}
s+ (1-\bar\varepsilon)\,\delta M^{2-p}\rho^p \leq t \leq s+ \delta M^{2-p}\rho^p,
\end{equation}
provided $K_{16\rho}(y) \times (s, s+ \delta M^{2-p}\rho^p] \subset E_T$.
\end{prop}
Again $\eta,\delta,\bar\varepsilon$ depend only on the data and $\alpha$, and they are not the constants of Proposition~\ref{exp pos}. In the singular case only the value at the center $y$ and the right endpoint of the time window are used below, so $\bar\varepsilon$ plays no role.

After the expansion of positivity, we recall Harnack Inequality. This requires the notation for the intrinsic cylinders as follows:  
$$(x_0, t_0) + Q_{\rho}^{+}(\theta) = K_{\rho}(x_0) \times (t_0, t_0 + \theta \rho^p],$$
$$(x_0, t_0) + Q_{\rho}^{-}(\theta) = K_{\rho}(x_0) \times (t_0 - \theta \rho^p, t_0].$$

\begin{thm}[Degenerate Harnack inequality]\label{thm:harnack}
Let $u$ be a continuous, nonnegative, weak solution to the
degenerate equation \eqref{p-Laplace}-\eqref{monotonicity} in $E_T$. There exist positive
constants $C_H$ and $\Gamma_H$, depending only on the data, such that for all 
intrinsic cylinders $(x_0, t_0) + Q_{4 \rho}^{\pm}(\theta) \subset E_T$ 
$$
\Gamma_H^{-1} \sup_{K_{\rho}(x_0)} u(\cdot, t_{0} - \theta \rho^p) \leq u(x_0, t_0) \leq \Gamma_H \inf_{K_{\rho}(x_0)} u(\cdot, t_{0} + \theta \rho^p),
$$
provided $\theta = \left( \frac{C_H}{u(x_0,t_0)} \right)^{p-2}.$
\end{thm}

\begin{thm}[Singular Harnack inequality]\label{thm:harnack sing}
Let $\frac{2N}{N+1}<p<2$, and let $u$ be a continuous, nonnegative, weak solution to the singular equation \eqref{p-Laplace}-\eqref{monotonicity}in $E_T$. There exist positive constants $C_H$ and $\Gamma_H$, depending only on the data, such that for all intrinsic cylinders $(x_0, t_0) + Q_{8 \rho}^{\pm}(\theta) \subset E_T$
$$
\Gamma_H^{-1} \sup_{K_{\rho}(x_0)} u(\cdot, t_{0} - \theta \rho^p) \leq u(x_0, t_0) \leq \Gamma_H \inf_{K_{\rho}(x_0)} u(\cdot, t_{0} + \theta \rho^p),
$$
provided $\theta = \left( \frac{C_H}{u(x_0,t_0)} \right)^{p-2}.$
\end{thm}
This statement requires the assumption $p>\frac{2N}{N+1}$, being $\frac{2N}{N+1}$ a \emph{critical} value. Indeed, when $p$ belongs to the subcritical range, namely $1<p \leq \frac{2N}{N+1}$, the Harnack inequality in this form is not available, and local boundedness of the solutions is not guaranteed either (see \cite{DGV}). Note that for $p<2$ the time lag $\theta\rho^p=(u(x_0,t_0)/C_H)^{2-p}\rho^p$ \emph{increases} with $u(x_0,t_0)$.

\begin{rem}\label{Rmk-Hoelder}
The Harnack inequality allows to prove the local H\"older continuity of the solutions (see \cite[Chapter~5, Section~4]{DGV}); the same holds in the singular supercritical range (see \cite{DGV}). Therefore, without loss of generality we may assume that the solutions we deal with are locally H\"older continuous, and therefore, their pointwise values are well-defined.
\end{rem}

We then state the Comparison Principle in bounded cylinders.   
\begin{prop}[Comparison Principle]
Let $E$ be a bounded open set in $\R^N$ with Lipschitz boundary $\partial E$ (a cube or a ball suffices). For $T > 0$, let $\mathbf{A}$ satisfy the conditions \eqref{structure} and \eqref{monotonicity} and let $u$, $v$ two solutions of \eqref{p-Laplace}-\eqref{monotonicity}such that
\begin{equation}
\begin{cases}
    u, v \in C (0, T; L^2(E)) \cap L^p (0, T; W^{1,p}(E)) \cap C (\bar{E}_T), \\
    u \geq v \text{ on } \partial_p E_{T}.
\end{cases}
\end{equation}
Then $u \geq v$ in $\bar{E}_T$.
\end{prop}
The proof relies on the monotonicity, as required in \eqref{monotonicity}, and holds both if $p>2$ and if $\frac{2N}{N+1}<p<2$. For a more general proof, see \cite[Chapter~7, Section~1]{DGV}. In conclusion, given that $\mathbf{A}$ is independent of $u$ and odd in the gradient variable, we have the following:
\begin{rem}\label{rem:shift}
If $u$ is a local weak solution to \eqref{p-Laplace}-\eqref{monotonicity}, then so are $u+\Theta$ and $-u+\Theta$ for every $\Theta\in\mathbb R$.
\end{rem}

The last two remarks of this section are concerned with an explicit solutions, and a scaling argument. It is relevant to note that these particular examples below hold both for the degenerate and the singular equation. 
\begin{rem}
\label{funzione paraboloidica}
For every $p>1$, 
\[
u(x,t)=t+ \frac{p-1}{p} \left( \frac{1}{N} \right)^{\frac{1}{p-1}}|x|^\frac{p}{p-1}
\]
is a nontrivial {\it eternal } solution to \eqref{Eq:prot}, i.e a solution in $\R^N \times \R$. With $A=\frac{p-1}{p}\big(\frac1N\big)^{\frac1{p-1}}$ one has $u<0$ in $\{t+A|x|^{\frac p{p-1}}<0\}$, alternative (i) of Theorem~\ref{P-L Thm in space} (Section~\ref{section3}) fails, and $u$ grows exactly like $|x|^{\frac p{p-1}}$. For $p<2$ the same function appears in Section~\ref{section7}, as $w_c$ with $c=\frac{p-1}pN^{-\frac1{p-1}}$ (so that $\kappa_c=1$). 
This particular function shows that the hypotheses below are non-vacuous and the growth rate is sharp.  
\end{rem}

\begin{rem}
\label{rescalings}
Let $a_1,a_2,a_3$ be real constants, and let $u$ be an eternal solution to \eqref{Eq:prot}. It is a matter of direct computations to show that for every $\lambda \in \R_+$, the rescaled function 
\[
v(x,t)=\lambda^{-a_1}u(\lambda^{a_2} x, \lambda^{a_3} t)
\] 
is still a solution to \eqref{Eq:prot} if $a_3=a_2p-a_1(p-2)$.
In addition, for every $(x,t) \in \mathcal{P}\coloneqq \{ t+A|x|^\frac{p}{p-1}<0\}$,  $(\lambda^{a_2} x, \lambda^{a_3} t) \in \mathcal{P}$ if $a_3=a_2 \frac{p}{p-1}$. 
As a consequence, if $u$ satisfies the hypotheses of Theorem~\ref{P-L Thm in space} with $h$ and $\Theta=0$, then $v$ satisfies them with $\lambda^{-a_3}h$ in place of $h$, provided that $a_2 \frac{p}{p-1}=a_3=a_2p-a_1(p-2)$.
\end{rem}

\section{A Theorem of Phragm\'en-Lindel\"of-Type in space} 
\label{section3}
In the following sections, we consider a truncated paraboloid defined as follows. Given some constants $A>0$, $h \ge 0$, let 
\begin{equation}
\label{Pa}
\mathcal{P}_h \coloneqq \set{(x,t) \in \R^N \times \R : t+A|x|^{\frac{p}{p-1}}<0, \, t<-h}. 
\end{equation}
\eqref{Pa} is a generalization of the set $\{(x,t) \in \R^N \times \R : t+A|x|^{2}<0, \, t<-h\}$ from \cite{Gla}. So, in such a similar setting, the aim of the following Theorem is to prove a certain growth rate with respect to the space variable, assuming $u$ negative in some $\Pa_h$. Actually, due to Remark~\ref{rem:shift}, this assumption is not far from requiring $u \le \Theta$ for some $\Theta \in \R$.

\begin{thm}
\label{P-L Thm in space}
Let $u$ be a local weak solution to the homogeneous degenerate equation \eqref{p-Laplace}-\eqref{monotonicity} in $\R^{N} \times (-\infty, 0]$. For some constants $A>0,$ $\Theta\in\R,$ $h \ge 0$, suppose that $u \le \Theta$ in $\mathcal{P}_h$. Then, the following alternative holds:
\begin{enumerate}
  \item[(i)] either $u(x, t) \le \Theta$ in the half-space $t<-h$,
  \item[(ii)] or 
  \begin{equation}
  \label{eq:Barenblatt_profile}
    \limsup_{\rho \to \infty} \left( \sup_{|x|=\rho, \, t<0} \frac{u(x,t)}{\rho^{\frac{p}{p-1}}} \right) > 0.
  \end{equation}
\end{enumerate}
\end{thm}
\begin{figure}
\vfill 
    \centering 
\begin{tikzpicture}[scale=1]
  \def\A{1}        
  \def\hh{0.8}     
  \def\pexp{1.5}   
  \pgfmathsetmacro{\xzero}{(\hh/\A)^(1/\pexp)}  
  \def\xmax{2.2}
  \pgfmathsetmacro{\tmin}{-\A*\xmax^\pexp}

  \draw[-{Latex[length=2mm]}] (-2.6,0) -- (2.6,0) node[right] {$x$};
  \draw[-{Latex[length=2mm]}] (0,\tmin-0.3) -- (0,0.5) node[above] {$t$};
  \node[below left] at (0,0) {\footnotesize $0$};

  \draw[dashed,gray] (-2.6,-\hh) -- (2.6,-\hh);
  \node[left] at (-2.6,-\hh) {\footnotesize $t=-h$};

  \fill[blue!12]
    (-\xmax,\tmin) --
    plot[domain=-\xmax:-\xzero,samples=60] ({\x},{-\A*(-\x)^\pexp}) --
    (-\xzero,-\hh) -- (\xzero,-\hh) --
    plot[domain=\xzero:\xmax,samples=60] ({\x},{-\A*(\x)^\pexp}) --
    (\xmax,\tmin) -- cycle;

  \draw[thick] (-\xzero,-\hh) -- (\xzero,-\hh);
  \draw[thick,domain=\xzero:\xmax,samples=60,smooth]
    plot ({\x},{-\A*(\x)^\pexp});
  \draw[thick,domain=-\xmax:-\xzero,samples=60,smooth]
    plot ({\x},{-\A*(-\x)^\pexp});

  \fill (\xzero,-\hh) circle (0.8pt);
  \fill (-\xzero,-\hh) circle (0.8pt);
  
  \node[right] at (1.55,-2.0) {\footnotesize $t=-A|x|^{\frac{p}{p-1}}$};
  \node[blue!70!black] at (0,-2.3)
    {\footnotesize $t+A|x|^{\frac{p}{p-1}}<0,\ t<-h$};

\end{tikzpicture}
\caption{An example of a truncated paraboloid $\mathcal{P}_h$.}
\label{fig:truncated paraboloid}
\end{figure}
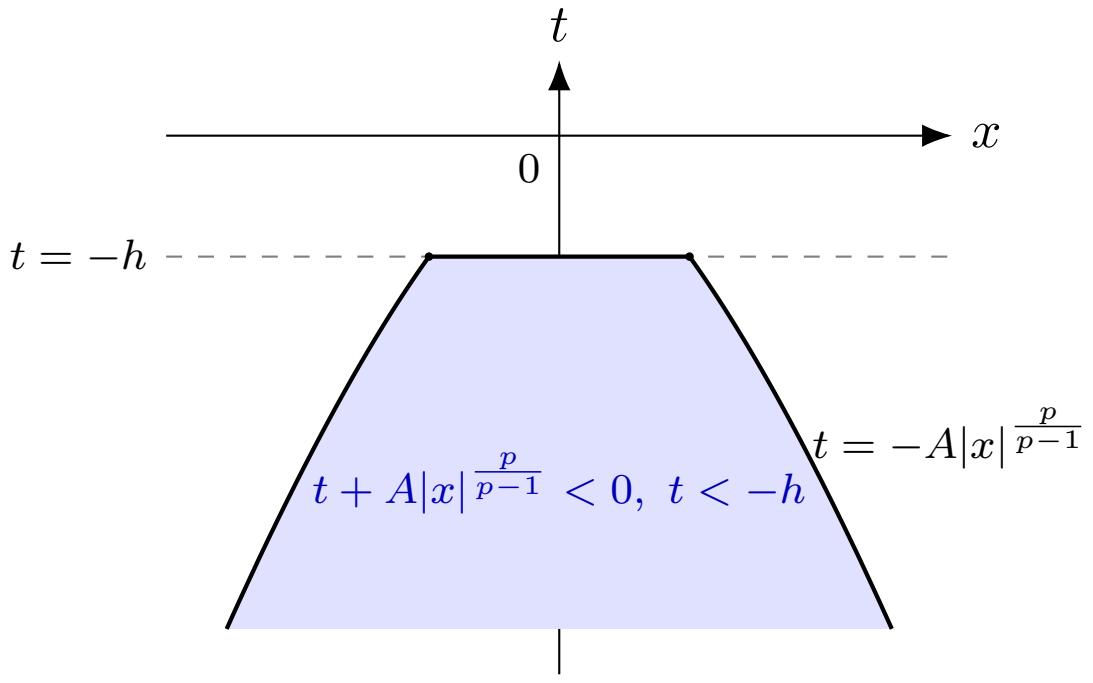

\begin{proof}
Without loss of generality, we may assume $\Theta=0$. By Remark~\ref{Rmk-Hoelder} we take $u$ to be continuous. Suppose that alternative (i) fails. Hence, there exists $(x_0, t_0) \notin \mathcal{P}_h$, with $t_0<-h$, such that $u_{0} \coloneqq u(x_0, t_0) > 0$. Let $\eta,\,\delta \in (0,1)$ and $b>1$ be the structural constants from Proposition~\ref{exp pos} with $\alpha=\frac12$ (the datum holds with $\alpha=1$, hence with $\frac12$), and set
\begin{align*}
\xi_\eta:=&\min\Big\{1,\tfrac{\eta}{1-\eta}\Big\},\qquad
\omega:=b^{p-2}\delta\,(2\eta)^{2-p}\,8^{-p},\\
C_{p,\xi_\eta}:=&\Big(1-(1+\xi_\eta)^{-\frac{p-1}{p}}\Big)^{-1},\qquad
\Lambda_\eta:=1+\sqrt N\,C_{p,\xi_\eta}.
\end{align*}
Here $\xi_\eta$ is just a real number, and it has nothing to do with the gradient variable $\xi\in\R^N$ of Section~\ref{sec:equations}. Fix
\begin{equation}\label{eq:rho0}
\rho_0\ \ge\ \max\left\{|x_0|,\ \left(\frac{A\Lambda_\eta^{\frac{p}{p-1}}}{\omega}\,u_0^{p-2}\right)^{\frac{p-1}{p(p-2)}}\right\}.
\end{equation}
We construct inductively $(x_k,t_k)$ with $u_k:=u(x_k,t_k)>0$, $t_{k+1}\le t_k$, and
\begin{equation}\label{eq:growth}
u_{k+1}\ge(1+\xi_\eta)u_k,\qquad \rho_k\le|x_{k+1}-x_k|\le\sqrt N\rho_k,\qquad
\rho_k:=\rho_0\Big(\tfrac{u_k}{u_0}\Big)^{\frac{p-1}{p}} .
\end{equation}
Given $(x_k,t_k)$, let $s_k:=t_k-\omega\rho_k^p u_k^{2-p}$ and
\[
Z^k:=K_{\rho_k}(x_k)\times(s_k,t_k],\qquad M_k:=\max_{\overline{Z^k}}u .
\]
\smallskip\noindent\emph{Claim~1: $\overline{K_{\rho_k}(x_k)}\times\{s_k\}\subset {\mathcal P}_h$.}

By \eqref{eq:growth}, $u_j\le(1+\xi_\eta)^{-(k-j)}u_k$, hence $\rho_j\le(1+\xi_\eta)^{-(k-j)\frac{p-1}{p}}\rho_k$ and, using $|x_0|\le\rho_0\le\rho_k$, for every $x\in K_{\rho_k}(x_k)$
\begin{equation}\label{eq:xbound}
|x|\le|x_0|+\sqrt N\sum_{j=0}^{k}\rho_j\le\Lambda_\eta\rho_k .
\end{equation}
Since $t_k<0$ we have $-s_k>\omega\rho_k^pu_k^{2-p}$, and $s_k<t_k\le t_0<-h$; hence it suffices that
\[
A(\Lambda_\eta\rho_k)^{\frac{p}{p-1}}\le\omega\rho_k^pu_k^{2-p},
\] 
i.e.
\[
A\Lambda_\eta^{\frac p{p-1}}\le\omega\,\rho_k^{\frac{p(p-2)}{p-1}}u_k^{2-p}.
\]
By \eqref{eq:growth}, $\rho_k^{\frac{p(p-2)}{p-1}}u_k^{2-p}=\rho_0^{\frac{p(p-2)}{p-1}}u_0^{2-p}$ is independent of $k$, and the inequality is exactly \eqref{eq:rho0}.

\smallskip\noindent\emph{Choice of the next vertex.}

By the Comparison Principle $M_k$ is attained on $\partial_pZ^k$. On the bottom $u\le0<u_k\le M_k$ by Claim~1, so that $M_k$ is necessarily attained at some $(x_{k+1},t_{k+1})$ on the lateral boundary: 
\[
|x_{k+1}-x_k|_\infty=\rho_k,\quad s_k<t_{k+1}\le t_k,\quad u_{k+1}=M_k.
\] 
This gives the second and third relations in \eqref{eq:growth}.

\smallskip\noindent\emph{Claim~2: $M_k\ge(1+\xi_\eta)u_k$.} 

We face a natural dichotomy: either $M_{k} \ge 2 u_k$ or $M_{k} < 2 u_k$. 

\smallskip\noindent\emph{Case 1: $M_k\ge2u_k$.} There is nothing to prove, since $\xi_\eta\le1$. 

\smallskip\noindent\emph{Case 2: $1\le M_k/u_k<2$.} Set
\[
r_k:=\frac{\rho_k}{8}\Big(\frac{M_k}{2u_k}\Big)^{\frac{p-2}{p}}<\frac{\rho_k}{8},
\]
so that $K_{8r_k}(x_k)\times(s_k,t_k]\subset Z^k$. The function $v:=M_k-u$ is a solution by Remark~\ref{rem:shift}, nonnegative in $Z^k$, and $v(\cdot,s_k)\ge M_k$ in $K_{r_k}(x_k)$ by Claim~1. Proposition~\ref{exp pos} with $\alpha=\frac12$, $M=M_k$, $\rho=r_k$ yields $v(\cdot,t)\ge\eta M_k$ in $K_{2r_k}(x_k)$ for
\[
s_k+\tfrac12\,b^{p-2}\delta(\eta M_k)^{2-p}r_k^p\ \le\ t\ \le\ s_k+b^{p-2}\delta(\eta M_k)^{2-p}r_k^p .
\]
A direct computation gives $b^{p-2}\delta(\eta M_k)^{2-p}r_k^p=\omega\rho_k^pu_k^{2-p}$, so the right endpoint is $t_k$. Hence $M_k-u_k=v(x_k,t_k)\ge\eta M_k$, i.e. $M_k\ge\frac{u_k}{1-\eta}\ge(1+\xi_\eta)u_k$.

\smallskip
We can now conclude. 
From \eqref{eq:growth}, 
\[
u_k\ge(1+\xi_\eta)^ku_0\to\infty,\,\,\text{ and }\,\,\rho_k\to\infty.
\] 
Suppose that a large constant $C$ bounds the whole sequence, that is, $|x_k|\le C$ for every $k$. Then $\rho_k\le|x_{k+1}-x_k|\le2C$, a contradiction; thus $|x_k|\to\infty$ along a subsequence. Evaluating \eqref{eq:xbound} at $x_k$ and inverting \eqref{eq:growth} yields
\begin{equation}\label{eq:final}
u_k=u_0\Big(\frac{\rho_k}{\rho_0}\Big)^{\frac p{p-1}}\ \ge\ \frac{u_0}{(\Lambda_\eta\rho_0)^{\frac p{p-1}}}\,|x_k|^{\frac p{p-1}},
\qquad t_k\le t_0<-h .
\end{equation}
Taking $\limsup$ along the subsequence gives alternative (ii), with the supremum restricted to $t<-h$.
\end{proof}

\begin{rem}
With $\xi_\eta$ fixed and $u_k\ge(1+\xi_\eta)^k u_0$ the construction we set up in the proof gives a quantitative lower bound $u_0(\Lambda_\eta\rho_0)^{-\frac p{p-1}}$ for the $\limsup$ in alternative (ii), depending only on the data, $A$, and $u_0$; it represents a sort of natural \emph{Phragm\'en-Lindel\"of} constant.
\end{rem}
\section{Weak profiles and self-similar solutions}
\label{section4}
After Theorem~\ref{P-L Thm in space}, one may look for a similar result with respect to the time variable. In this section we show that \eqref{Eq:prot} has \emph{self-similar} solutions
\[
u(x,t)=\tau^{\beta}F(x\tau^{-\gamma}),\qquad \tau:=-t,
\]
which are positive in $\Pa_h$ and behave like $T^\beta$ for every $\beta\le 0$ and every $\beta\ge1$ (for $\beta=1$ provided the amplitude is large enough, for $\beta>1$ provided $h$ is large), whereas none is positive for $0<\beta<1$. Hence positivity alone does not bound the growth in time from above, and a result in time is necessarily a dichotomy, the object of Section~\ref{section5}.

We now introduce the following notation. We set $\phi_p(s):=|s|^{p-2}s$ and $\psi(s):=\operatorname{sign}(s)\,|s|^{\frac1{p-1}}$, the inverse of $\phi_p$. The function $\psi$ is continuous, though not Lipschitz at $0$. Given $A>0$ and $h\ge0$, let
\begin{equation}\label{Time slice}
R(\tau):=\Big(\frac{\tau}{A}\Big)^{\frac{p-1}{p}},\qquad \tau>h,
\end{equation}
be the radius of the slice of $\Pa_h$ at time $-\tau$. For $\beta\in\R$ we put
\begin{equation}\label{gamma}
\gamma=\gamma_\beta:=\frac{1+\beta(p-2)}{p}.
\end{equation}
Two identities follow at once from \eqref{gamma} and are used below:
\begin{equation}\label{gamma-id}
(\beta-\gamma)(p-1)-\gamma=\beta-1,\qquad
\gamma-\frac{p-1}{p}=\frac{(\beta-1)(p-2)}{p}.
\end{equation}

\begin{defn}\label{def:profile}
Let $\beta\in\R$ and $\gamma$ as in \eqref{gamma}. A function $F\in C^1(\R^N)$ is a \emph{weak profile} if
\begin{equation}\label{weak profile}
\int_{\R^N}\Big[\phi_p( D F)\cdot D\varphi+\big(\gamma\,y\cdot D F-\beta F\big)\varphi\Big]\,dy=0
\qquad\text{for every }\varphi\in C^{1}_c(\R^N).
\end{equation}
\end{defn}

\eqref{weak profile} is the equation $\operatorname{div}\phi_p( D F)=\gamma\,y\cdot D F-\beta F$ in a weak sense. The next two theorems answer two questions: for which exponents $\beta$ a radial weak profile exists, and when the resulting solution $u$ is strictly positive in $\Pa_h$.

\begin{thm}\label{thm:profile}
For every $\beta\in\R$ there exists $f\in C^1([0,\infty))$ with $f(0)=1$ and $f'(0)=0$ such that $g:=\phi_p(f')$ is continuous, $g(0)=0$, $|g(r)|\le Cr$ near $0$, and
\begin{equation}\label{ode}
f'=\psi(g),\qquad (r^{N-1}g)'=r^{N-1}\big(\gamma\,r\,\psi(g)-\beta f\big)\quad\text{in }(0,\infty).
\end{equation}
Moreover $F(y):=f(|y|)$ is a radial weak profile. Let $y_1:=\sup\{r>0:\ f>0\text{ on }[0,r)\}$. Then:
\begin{enumerate}
\item[(a)] if $\beta<0$, then $g>0$ in $(0,\infty)$, so $f$ is strictly increasing, $f\ge1$ and $y_1=\infty$;
\item[(b)] if $\beta=0$, we take $f\equiv1$, and $y_1=\infty$;
\item[(c)] if $\beta>0$, then $g<0$ in $(0,y_1]$, $f$ is strictly decreasing in $[0,y_1]$, $y_1<\infty$ and $f'(y_1)<0$. Moreover, with $\Lambda_\beta:=2\beta/N$,
\[
y_1\ \ge\ \min\left\{\rho_\beta,\ \left(\frac{p}{(p-1)\Lambda_\beta^{\frac1{p-1}}}\right)^{\frac{p-1}{p}}\right\},\qquad
\rho_\beta:=\left(\frac{\beta\,(N+\frac{p}{p-1})}{2N\gamma\,\Lambda_\beta^{\frac1{p-1}}}\right)^{\frac{p-1}{p}};
\]
\item[(d)] if $\beta=1$, then $f(r)=1-k\,r^{\frac{p}{p-1}}$ with $k=\frac{p-1}{p}N^{-\frac1{p-1}}$, and $y_1=k^{-\frac{p-1}{p}}$;
\item[(e)] if $\beta\neq0$, every continuous solution of the integral system \eqref{integral} below coincides with $(f,g)$ on the common interval of definition, intersected with $[0,y_1]$ (with $y_1=\infty$ if $\beta<0$).
\end{enumerate}
For $\beta=0$ the integral system \eqref{integral} also has nonconstant solutions, namely $g=C_\sharp\, r^{\frac{2(p-1)}{p-2}}$ with $C_\sharp^{\frac{p-2}{p-1}}=\gamma\big/\big(N+1+\frac2{p-2}\big)$, $\gamma=\frac1p$; this is why (e) requires $\beta\ne0$, and for $\beta=0$ we choose $f\equiv1$.
\end{thm}

\begin{thm}\label{thm:positivity}
Let $f=f_\beta$ be as in Theorem~\ref{thm:profile}. For $c>0$ and $s\ge0$ let $\nu_c:=c^{-\frac{p-2}{p}}$ and, for $(x,t)\in\R^N\times(-\infty,s)$ and $\tau:=-t$,
\[
u_{\beta,c,s}(x,t):=c\,(\tau+s)^{\beta}\,f\big(\nu_c\,|x|\,(\tau+s)^{-\gamma}\big).
\]
Then $u_{\beta,c,s}$ is a weak solution of \eqref{Eq:prot} in $\R^N\times(-\infty,s)$ and $u_{\beta,c,s}(0,-T)=c\,(T+s)^{\beta}$. Moreover, let $y_c:=c^{\frac{p-2}{p}}y_1\in(0,\infty]$ and
\[
\Phi_s(\tau):=R(\tau)\,(\tau+s)^{-\gamma}.
\]
Then $u_{\beta,c,s}>0$ in $\Pa_h$ \emph{if and only if} $\sup_{\tau>h}\Phi_s(\tau)\le y_c$. Consequently, for every $A>0$, $h\ge0$, $c>0$, $s\ge0$:
\begin{enumerate}
\item[(1)] if $\beta\le0$, then $u_{\beta,c,s}$ is positive in all of $\R^N\times(-\infty,s)$.
\item[(2)] if $0<\beta<1$, then $u_{\beta,c,s}$ is \emph{never} positive in $\Pa_h$.
\item[(3)] if $\beta=1$, then $u_{1,c,s}$ is positive in $\Pa_h$ if and only if $c\ge a(A)$, where
\[
a(A):=\Big(N\Big(\frac{p}{p-1}\Big)^{p-1}\Big)^{-\frac1{p-2}}A^{-\frac{p-1}{p-2}} .
\]
\item[(4)] if $\beta>1$, then $u_{\beta,c,0}$ is positive in $\Pa_h$ if and only if $h>0$ and $h\ge h_0(c):=h_0(1)\,c^{-\frac1{\beta-1}}$, where $h_0(1)=\big(A^{-\frac{p-1}{p}}/y_1\big)^{\frac{p}{(\beta-1)(p-2)}}$. For $s\ge0$, the condition $h+s\ge h_0(c)$ is sufficient.
\end{enumerate}
\end{thm}

\begin{rem}\label{rem:beta1}
By the proof of Theorem~\ref{thm:positivity}, $u_{1,c,0}(x,t)=c\,\tau-c^{\frac1{p-1}}k\,|x|^{\frac{p}{p-1}}$, which is, up to the change $u\mapsto-u$ of Remark~\ref{rem:shift} and a scaling, the function of Remark~\ref{funzione paraboloidica}. Moreover, for every $h>0$, $s\ge0$ and $\beta>1$, the choice $c\ge\big(h_0(1)/(h+s)\big)^{\beta-1}$ gives a solution positive in $\Pa_h$ with $u(0,-T)=c\,(T+s)^{\beta}$.
\end{rem}

\subsection{Reduction to a profile equation}
\begin{lem}\label{lem:ansatz}
Let $\beta\in\R$, $F\in C^1(\R^N)$, and $u(x,t):=\tau^{\beta}F(x\tau^{-\gamma})$ for $t<0$. Then $u\in C^1(\R^N\times(-\infty,0))$ and, with $y=x\tau^{-\gamma}$,
\begin{equation}\label{ansatz}
u_t=\tau^{\beta-1}\big(\gamma\,y\cdot D F-\beta F\big)(y),\qquad
\phi_p( D u)=\tau^{(\beta-\gamma)(p-1)}\,\phi_p( D F)(y).
\end{equation}
If $F$ is a weak profile, then $u$ is a weak solution of \eqref{Eq:prot} in $\R^N\times(-\infty,0)$.
\end{lem}
\begin{proof}
\emph{Formulas.} Since $\partial_\tau\big[F(x\tau^{-\gamma})\big]=-\gamma\tau^{-1}\,y\cdot D F(y)$, we get $u_\tau=\tau^{\beta-1}(\beta F-\gamma\,y\cdot D F)$, and $u_t=-u_\tau$. Next $ D u=\tau^{\beta-\gamma} D F(y)$, and $\phi_p$ is homogeneous of degree $p-1$.

\emph{Weak formulation.} Fix $t<0$ and $\varphi\in C^1_c(\R^N)$, and put $\tilde\varphi(y):=\varphi(\tau^{\gamma}y)$. Changing variables $x=\tau^\gamma y$, so that $ D_y\tilde\varphi=\tau^{\gamma} D\varphi$ and $dx=\tau^{N\gamma}dy$, and using the first part,
\[
\int_{\R^N}\phi_p( D u)\cdot D\varphi\,dx
=\tau^{(\beta-\gamma)(p-1)-\gamma+N\gamma}\int_{\R^N}\phi_p( D F)\cdot D\tilde\varphi\,dy .
\]
By \eqref{weak profile} the last integral equals $-\int_{\R^N}(\gamma\,y\cdot D F-\beta F)\tilde\varphi\,dy$. Going back to the variable $x$, and using \eqref{gamma-id} to see that $(\beta-\gamma)(p-1)-\gamma=\beta-1$, we obtain
\begin{equation}\label{flux identity}
\int_{\R^N}\phi_p( D u)\cdot D\varphi\,dx=-\int_{\R^N}u_t\,\varphi\,dx\qquad\text{for every }t<0 .
\end{equation}
Since $u\in C^1$, for $[t_1,t_2]\subset(-\infty,0)$ and for every test function $\varphi$ as in Definition~\ref{def:weak} (after a time translation, and with $E=\R^N$) we have $\int_{\R^N}u\varphi\,dx\big|_{t_1}^{t_2}=\iint(u_t\varphi+u\varphi_t)\,dx\,dt$, so the integral identity in the definition of weak solution reduces to $\iint(u_t\varphi+\phi_p( D u)\cdot D\varphi)\,dx\,dt$, which vanishes by \eqref{flux identity}. Both $u_t$ and $\phi_p( D u)$ are locally bounded, so \eqref{flux identity} extends by density from $C^1_c(\R^N)$ to the functions $\varphi(\cdot,t)\in W^{1,p}(\R^N)$ with compact support, and the same holds for the test functions in the definition of weak solution.
\end{proof}

\begin{lem}\label{lem:sym}
Let $F$ be a weak profile. Then, for every $\lambda>0$, $F_\lambda(y):=\lambda F(\nu_\lambda y)$, where $\nu_\lambda:=\lambda^{-\frac{p-2}{p}}$, is a weak profile. If $F$ is radial with first zero $y_1$, then $F_\lambda$ has first zero $\lambda^{\frac{p-2}{p}}y_1$. Moreover, if $u$ is a weak solution of \eqref{Eq:prot}, so is $u(x,t-s)$ for every $s\in\R$.
\end{lem}
\begin{proof}
Note that $\lambda^{p-1}\nu_\lambda^{p}=\lambda$. Let $\varphi\in C^1_c(\R^N)$ and put $\tilde\varphi(z):=\varphi(z/\nu_\lambda)$, so that $ D\varphi(y)=\nu_\lambda\, D\tilde\varphi(\nu_\lambda y)$. We have $ D F_\lambda(y)=\lambda\nu_\lambda\, D F(\nu_\lambda y)$ and $y\cdot D F_\lambda(y)=\lambda\,z\cdot D F(z)$ with $z=\nu_\lambda y$. The change of variables $z=\nu_\lambda y$, $dy=\nu_\lambda^{-N}dz$, gives
\[
\begin{gathered}
\int_{\R^N}\phi_p( D F_\lambda)\cdot D\varphi\,dy=\lambda^{p-1}\nu_\lambda^{p}\,\nu_\lambda^{-N}I,\\
\int_{\R^N}\big(\gamma\,y\cdot D F_\lambda-\beta F_\lambda\big)\varphi\,dy=\lambda\,\nu_\lambda^{-N}J,
\end{gathered}
\]
where $I:=\int\phi_p( D F)\cdot D\tilde\varphi\,dz$ and $J:=\int(\gamma\,z\cdot D F-\beta F)\tilde\varphi\,dz$. By \eqref{weak profile} for $F$ and the test function $\tilde\varphi$, $I+J=0$. Since $\lambda^{p-1}\nu_\lambda^p=\lambda$, the left-hand side of \eqref{weak profile} for $F_\lambda$ is $\lambda\nu_\lambda^{-N}(I+J)=0$. The statement on the zero is immediate, and the invariance under time translations is clear.
\end{proof}

\subsection{Theorem~\ref{thm:profile}: existence of the profile}
For $F(y)=f(|y|)$ and $g=\phi_p(f')$, the profile equation becomes \eqref{ode}. We solve it first near $r=0$, as the fixed point of an integral system, then continue the solution to $[0,\infty)$.

\begin{lem}\label{lem:local}
Let $\beta\in\R$ and $M_\beta:=\big(\tfrac32|\beta|+1\big)/N$. There exist $r_0>0$ and $f,g\in C([0,r_0])$ with $|g(r)|\le M_\beta r$ and $f\in[\frac12,\frac32]$, such that
\begin{equation}\label{integral}
g(r)=r^{1-N}\!\int_0^r s^{N-1}\big(\gamma\,s\,\psi(g(s))-\beta f(s)\big)\,ds,\qquad f(r)=1+\int_0^r\psi(g(s))\,ds .
\end{equation}
Moreover $f\in C^1$, $f'(0)=0$, $r^{N-1}g\in C^1$, and \eqref{ode} holds in $(0,r_0]$.
\end{lem}
\begin{proof}
Choose $r_0>0$ such that
\[
(\mathrm R_1)\ \ \frac{p-1}{p}\,M_\beta^{\frac1{p-1}}r_0^{\frac{p}{p-1}}\le\frac12,\qquad
(\mathrm R_2)\ \ |\gamma|\,r_0\,(M_\beta r_0)^{\frac1{p-1}}\le1 .
\]
Let $K_\beta:=\{g\in C([0,r_0]):\ |g(r)|\le M_\beta r\}$, a closed convex subset of $C([0,r_0])$, with the sup norm. For $g\in K_\beta$ put
\[
\begin{gathered}
f_g(r):=1+\int_0^r\psi(g),\qquad \ell_g(s):=\gamma\,s\,\psi(g(s))-\beta f_g(s),\\
(\mathcal Sg)(r):=r^{1-N}\int_0^rs^{N-1}\ell_g(s)\,ds .
\end{gathered}
\]
\begin{itemize}
\item \emph{$\mathcal S(K_\beta)\subset K_\beta$.} Since $|\psi(g(s))|\le(M_\beta s)^{\frac1{p-1}}$, we have 
\[
|f_g-1|\le\frac{p-1}{p}M_\beta^{\frac1{p-1}}r^{\frac{p}{p-1}}\le\frac12
\] 
by $(\mathrm R_1)$. By $(\mathrm R_2)$, $|\ell_g|\le1+\frac32|\beta|=NM_\beta$. So $|\mathcal Sg(r)|\le r\,\sup|\ell_g|/N\le M_\beta r$.
\item \emph{$\mathcal S(K_\beta)$ is relatively compact.} For $r>0$, $(\mathcal Sg)'=(1-N)r^{-N}\int_0^rs^{N-1}\ell_g+\ell_g$ is bounded by $(2N-1)M_\beta$ uniformly in $g\in K_\beta$. So $\mathcal S(K_\beta)$ is equi-Lipschitz, and Ascoli--Arzel\`a applies.
\item \emph{$\mathcal S$ is continuous.} If $g_n\to g$ uniformly, then $\psi(g_n)\to\psi(g)$ uniformly, because $\psi$ is uniformly continuous on bounded sets. Hence $f_{g_n}\to f_g$ and $\ell_{g_n}\to\ell_g$ uniformly, and so $\mathcal Sg_n\to\mathcal Sg$.
\end{itemize}
By Schauder's fixed point theorem, $\mathcal S$ has a fixed point $g\in K_\beta$; set $f:=f_g$. The regularity follows from \eqref{integral}: $f'=\psi(g)$, $f'(0)=0$, and $r^{N-1}g=\int_0^rs^{N-1}\ell_g$ is $C^1$ with derivative $r^{N-1}\ell_g(r)$.
\end{proof}

\begin{lem}\label{lem:origin}
Every continuous solution $(f,g)$ of \eqref{integral} on some $[0,r_1]$ satisfies
\[
g(r)=-\frac{\beta}{N}\,r+O\big(r^{1+\frac{p}{p-1}}\big),\qquad f(r)=1+O\big(r^{\frac{p}{p-1}}\big)\qquad(r\to0).
\]
\end{lem}
\begin{proof}
The integrand in \eqref{integral} is bounded near $0$, hence $|g(r)|\le C r$. Then $|\psi(g)|\le Cr^{\frac1{p-1}}$, so $f=1+O(r^{\frac{p}{p-1}})$ and $\gamma\,s\,\psi(g(s))=O(s^{\frac{p}{p-1}})$. Inserting this into \eqref{integral},
\[
g(r)=r^{1-N}\int_0^rs^{N-1}\big(-\beta+O(s^{\frac{p}{p-1}})\big)\,ds=-\frac\beta N\,r+O\big(r^{1+\frac{p}{p-1}}\big).\qedhere
\]
\end{proof}

\begin{lem}\label{lem:global}
The solution of Lemma~\ref{lem:local} extends to a solution of \eqref{ode} in $(0,\infty)$.
\end{lem}
\begin{proof}
For $r\ge r_0$ consider the first order system
\begin{equation}\label{system}
f'=\psi(g),\qquad g'=-\frac{N-1}{r}\,g+\gamma\,r\,\psi(g)-\beta f ,
\end{equation}
with initial datum $(f(r_0),g(r_0))$. Its right-hand side is continuous on $[r_0,\infty)\times\R^2$, so by Peano's theorem there is a maximal solution on $[r_0,r_{\max})$. If $r_{\max}<\infty$, the solution leaves every compact subset of $[r_0,\infty)\times\R^2$ as $r\to r_{\max}$, that is
\begin{equation}\label{blowup}
\limsup_{r\to r_{\max}}\big(|f(r)|+|g(r)|\big)=\infty .
\end{equation}
We show that \eqref{blowup} is impossible. Let $w:=|f|+|g|$, which is continuous on $[r_0,r_{\max})$. For $r\in[r_0,r_{\max})$,
\[
\begin{gathered}
|f(r)|\le|f(r_0)|+\int_{r_0}^r|\psi(g)|\,ds,\\
|g(r)|\le|g(r_0)|+\int_{r_0}^r\Big(\frac{N-1}{s}|g|+|\gamma|\,s\,|\psi(g)|+|\beta|\,|f|\Big)ds .
\end{gathered}
\]
Since $\frac1{p-1}<1$ we have $|\psi(z)|=|z|^{\frac1{p-1}}\le1+|z|$ for every $z\in\R$; moreover $\frac1s\le\frac1{r_0}$ and $s\le r_{\max}$ on the interval. Adding the two inequalities,
\[
\begin{gathered}
w(r)\ \le\ w(r_0)+L_1\,(r_{\max}-r_0)+L_2\int_{r_0}^rw(s)\,ds,\\
L_1:=1+|\gamma|\,r_{\max},\quad L_2:=\max\Big\{1+\frac{N-1}{r_0}+|\gamma|\,r_{\max},\,|\beta|\Big\}.
\end{gathered}
\]
By Gronwall's lemma in integral form
\[
w(r)\ \le\ \big(w(r_0)+L_1(r_{\max}-r_0)\big)\,e^{L_2(r_{\max}-r_0)}\qquad\text{for all }r\in[r_0,r_{\max}),
\]
which contradicts \eqref{blowup}. Hence $r_{\max}=\infty$. Finally, for $r>r_0$ the system \eqref{system} is equivalent to \eqref{ode}, because $(r^{N-1}g)'=r^{N-1}\big(g'+\frac{N-1}{r}g\big)$; and the glued function satisfies \eqref{ode} also at $r_0$, since the one-sided derivatives of $r^{N-1}g$ at $r_0$ are both equal to $r_0^{N-1}\big(\gamma\,r_0\psi(g(r_0))-\beta f(r_0)\big)$.
\end{proof}

\begin{lem}\label{lem:wprofile}
$F(y):=f(|y|)$ is a weak profile.
\end{lem}
\begin{proof}
Since $f'(0)=0$, $F\in C^1(\R^N)$, with $ D F(y)=f'(|y|)\,y/|y|$ and $ D F(0)=0$. Let $G:=\phi_p( D F)=g(|y|)\,y/|y|$; it is continuous, because $|g(r)|\le Cr$ near $0$, and it is of class $C^1$ in $\R^N\setminus\{0\}$, because $g=r^{1-N}(r^{N-1}g)$ is $C^1$ for $r>0$. By \eqref{ode}, for $y\ne0$,
\[
\operatorname{div}G(y)=|y|^{1-N}\big(r^{N-1}g\big)'(|y|)=\gamma\,|y|\psi(g(|y|))-\beta f(|y|)=\gamma\,y\cdot D F-\beta F ,
\]
where we used $y\cdot D F=r f'(r)=r\psi(g)$. Let $\varphi\in C^1_c(\R^N)$ and $\varepsilon>0$. The divergence theorem on $\{|y|>\varepsilon\}$, whose outer normal on $\{|y|=\varepsilon\}$ is $-y/|y|$, gives
\[
\int_{\{|y|>\varepsilon\}}G\cdot D\varphi\,dy=-\int_{\{|y|>\varepsilon\}}\big(\gamma\,y\cdot D F-\beta F\big)\varphi\,dy-g(\varepsilon)\int_{\{|y|=\varepsilon\}}\varphi\,d\sigma .
\]
The last term is bounded by $C\varepsilon\cdot\varepsilon^{N-1}\|\varphi\|_\infty\to0$, and the other integrands are bounded near the origin. Letting $\varepsilon\to0$ gives \eqref{weak profile}.
\end{proof}
Lemmas~\ref{lem:local}, \ref{lem:global} and \ref{lem:wprofile} prove the existence part of Theorem~\ref{thm:profile}.

\subsection{Theorem~\ref{thm:profile}: qualitative properties}
We use repeatedly the following observation, which follows from \eqref{ode}: at a zero $r_g>0$ of $g$ one has $\psi(g(r_g))=0$, hence
\begin{equation}\label{zero of g}
(r^{N-1}g)'(r_g)=-\beta\,r_g^{N-1}f(r_g).
\end{equation}

\begin{prop}[the case $\beta\le0$]\label{prop:neg}
If $\beta<0$, then $g>0$ in $(0,\infty)$. Hence $f$ is strictly increasing, $f\ge1$ and $y_1=\infty$. If $\beta=0$, then $f\equiv1$ is a solution.
\end{prop}
\begin{proof}
By Lemma~\ref{lem:origin}, $g=\frac{|\beta|}{N}r+o(r)>0$ near $0$. Suppose $g$ has a first zero $r_g$. On $(0,r_g)$ we have $f'=\psi(g)>0$, so $f(r_g)\ge1$, and \eqref{zero of g} gives $(r^{N-1}g)'(r_g)=|\beta|r_g^{N-1}f(r_g)>0$. But $r^{N-1}g>0$ on $(0,r_g)$ and vanishes at $r_g$, which forces $(r^{N-1}g)'(r_g)\le0$: a contradiction. The case $\beta=0$ is immediate from \eqref{ode}.
\end{proof}

\begin{prop}[the case $\beta>0$]\label{prop:pos}
Let $\beta>0$. Then $g<0$ in $(0,y_1]$, $f$ is strictly decreasing in $[0,y_1]$, $y_1<\infty$, $f'(y_1)<0$, and the lower bound on $y_1$ in Theorem~\ref{thm:profile}(c) holds.
\end{prop}
\begin{proof}
Here $\gamma>0$. Write $G:=-g$ and $H:=r^{N-1}G$. By \eqref{ode},
\begin{equation}\label{eq for H}
H'=r^{N-1}\big(\gamma\,r\,G^{\frac1{p-1}}+\beta f\big)\quad\text{where }G\ge0 .
\end{equation}

\emph{Sign of $g$.} By Lemma~\ref{lem:origin}, $g<0$ near $0$. If $g$ had a first zero $r_g<y_1$, then $f(r_g)>0$ and \eqref{zero of g} would give $(r^{N-1}g)'(r_g)<0$, whereas $r^{N-1}g<0$ on $(0,r_g)$ and vanishes at $r_g$, so the derivative at $r_g$ is $\ge0$. Hence $g<0$ in $(0,y_1)$ and $f'=\psi(g)<0$ there.

\emph{$y_1<\infty$.} Suppose $f>0$ in $[0,\infty)$; then $G>0$ in $(0,\infty)$. Since $\beta f>0$, \eqref{eq for H} gives $H'>0$ and
\[
H'\ \ge\ \gamma\,r^{N-\frac{N-1}{p-1}}\,H^{\frac1{p-1}},\qquad\text{because }r^{N-1}\,r\,G^{\frac1{p-1}}=r^{N-\frac{N-1}{p-1}}H^{\frac1{p-1}}.
\]
Fix $r_*>0$, so that $H(r_*)>0$. Dividing by $H^{\frac1{p-1}}$ and integrating over $(r_*,r)$, with $1-\frac1{p-1}=\frac{p-2}{p-1}$ and $\kappa:=N+1-\frac{N-1}{p-1}$,
\[
\begin{aligned}
H(r)^{\frac{p-2}{p-1}}\ &\ge\ H(r_*)^{\frac{p-2}{p-1}}+\frac{p-2}{p-1}\,\gamma\,\frac{r^{\kappa}-r_*^{\kappa}}{\kappa}\\
&\ge\ C\,r^{\kappa}\qquad(r\ge2r_*).
\end{aligned}
\]
Raising to the power $\frac{p-1}{p-2}$ we get $H\ge C\,r^{N+\frac{p}{p-2}}$, hence $G=r^{1-N}H\ge C\,r^{1+\frac{p}{p-2}}$ and
\[
f'=-G^{\frac1{p-1}}\ \le\ -C\,r^{\frac{2}{p-2}} .
\]
Integrating, $f(r)\to-\infty$, a contradiction. So $y_1<\infty$.

\emph{Transversality.} By \eqref{eq for H}, $H'>0$ in $(0,y_1)$, and $H(0)=0$; hence $H(y_1)>0$, that is $g(y_1)<0$ and $f'(y_1)=\psi(g(y_1))<0$.

\emph{Lower bound for $y_1$.} In $(0,y_1)$ we have $0<f\le1$, so \eqref{eq for H} gives
\begin{equation}\label{bound G}
G(r)\ \le\ r^{1-N}\int_0^rs^{N-1}\big(\gamma\,s\,G(s)^{\frac1{p-1}}+\beta\big)\,ds .
\end{equation}
Let $\Lambda_\beta=2\beta/N$. Suppose $G(s)\le\Lambda_\beta s$ for $s\in[0,r]$ and $r\le\rho_\beta$. Then \eqref{bound G} yields
\[
G(r)\ \le\ r\Big[\frac\beta N+\frac{\gamma\Lambda_\beta^{\frac1{p-1}}\,r^{\frac{p}{p-1}}}{N+\frac{p}{p-1}}\Big]\ \le\ \frac{3\beta}{2N}\,r\ <\ \Lambda_\beta\,r,
\]
by the definition of $\rho_\beta$. Since $G(r)=\frac\beta Nr+o(r)<\Lambda_\beta r$ for small $r>0$, a continuity argument shows that $G(r)\le\Lambda_\beta r$ in $[0,\min\{\rho_\beta,y_1\})$: at a first point where $G=\Lambda_\beta r$ we would have just proved $G<\Lambda_\beta r$. Then
\[
f(r)=1-\int_0^rG^{\frac1{p-1}}\ \ge\ 1-\frac{p-1}{p}\,\Lambda_\beta^{\frac1{p-1}}\,r^{\frac{p}{p-1}},
\]
which is positive for $r<\big(\frac{p}{(p-1)\Lambda_\beta^{\frac1{p-1}}}\big)^{\frac{p-1}{p}}$. Hence $y_1$ is not smaller than the minimum of the two numbers.
\end{proof}

\begin{prop}[uniqueness]\label{prop:uniq}
Let $\beta\ne0$ and let $(f_2,g_2)$ be a continuous solution of \eqref{integral} on some $[0,r_1]$. Then $(f_2,g_2)=(f,g)$ on $[0,\min\{r_1,y_1\}]$, where $y_1=\infty$ if $\beta<0$.
\end{prop}
\begin{proof}
Write $(f_1,g_1)=(f,g)$.

\emph{Near the origin.} By Lemma~\ref{lem:origin} and $\beta\ne0$ there are $C,r_2>0$ such that $g_1,g_2$ have the same sign and $|g_i(r)|\ge r/C$ on $(0,r_2]$. For every $r\in(0,r_2]$, the function $\psi$ is differentiable between $g_1(r)$ and $g_2(r)$, so by the mean value theorem $\psi(g_1(r))-\psi(g_2(r))=\psi'(\xi_r)\,\big(g_1(r)-g_2(r)\big)$, where $\xi_r$ is some real number between $g_1(r)$ and $g_2(r)$. Since $|\xi_r|\ge r/C$, we get $\psi'(\xi_r)=\frac1{p-1}|\xi_r|^{-\frac{p-2}{p-1}}\le C\,r^{-\frac{p-2}{p-1}}$, hence
\[
|\psi(g_1)-\psi(g_2)|(r)\ \le\ C\,r^{-\frac{p-2}{p-1}}\,|g_1-g_2|(r).
\]
Let $E(r):=\max_{[0,r]}|g_1-g_2|$. The second equation in \eqref{integral} gives $|f_1-f_2|(r)\le C\,r^{\frac1{p-1}}E(r)$, because $\int_0^rs^{-\frac{p-2}{p-1}}ds=(p-1)r^{\frac1{p-1}}$. The first equation then gives
\[
\begin{aligned}
|g_1-g_2|(r)&\le r^{1-N}\!\int_0^rs^{N-1}\Big(|\gamma|\,s\,|\psi(g_1)-\psi(g_2)|+|\beta|\,|f_1-f_2|\Big)ds\\
&\le\ C\,r^{\frac{p}{p-1}}E(r),
\end{aligned}
\]
where we used $s\cdot s^{-\frac{p-2}{p-1}}=s^{\frac1{p-1}}$. The right-hand side is nondecreasing in $r$, so $E(r)\le C\,r^{\frac{p}{p-1}}E(r)$, and $E\equiv0$ on a small interval $[0,r_3]$.

\emph{Away from the origin.} Let $[0,\bar r]$ be the maximal interval on which the two solutions coincide, and suppose $\bar r<\min\{r_1,y_1\}$. By Propositions~\ref{prop:neg} and~\ref{prop:pos}, $g(\bar r)\ne0$. Both solutions solve \eqref{system} on $[r_3,r_1]$, and the right-hand side of \eqref{system} is locally Lipschitz near $\big(f(\bar r),g(\bar r)\big)$, where $\psi$ is smooth. Uniqueness for the Cauchy problem at $\bar r$ contradicts the maximality of $\bar r$.
\end{proof}

\begin{prop}[the case $\beta=1$]\label{prop:beta1}
For $\beta=1$ one has $f(r)=1-k\,r^{\frac{p}{p-1}}$ with $k=\frac{p-1}{p}N^{-\frac1{p-1}}$, and $y_1=k^{-\frac{p-1}{p}}$.
\end{prop}
\begin{proof}
Let $f=1-k\,r^{\frac{p}{p-1}}$. Then $f'=-\frac{p}{p-1}k\,r^{\frac1{p-1}}$ and $\phi_p(f')=-\big(\tfrac{p}{p-1}k\big)^{p-1}r=-\frac rN$, by the choice of $k$. So $g=-r/N$ and $(r^{N-1}g)'=-r^{N-1}$. With $\beta=1$ and $\gamma=\frac{p-1}{p}$, the right-hand side of \eqref{ode} is
\[
r^{N-1}\Big(\frac{p-1}{p}\,r\cdot\big(-\tfrac{p}{p-1}k\,r^{\frac1{p-1}}\big)-1+k\,r^{\frac{p}{p-1}}\Big)=-r^{N-1}.
\]
So $(f,g)$ solves \eqref{ode}, with $f(0)=1$ and $r^{N-1}g\to0$ as $r\to0$; integrating \eqref{ode} on $(0,r)$ shows that it solves \eqref{integral}. By Proposition~\ref{prop:uniq} it is the profile. Its first zero is given by $k\,r^{\frac{p}{p-1}}=1$.
\end{proof}

Propositions~\ref{prop:neg}--\ref{prop:beta1} prove parts (a)--(e) of Theorem~\ref{thm:profile}.

\subsection{Proof of Theorem~\ref{thm:positivity}}
\emph{Weak solution.} By Lemma~\ref{lem:wprofile} and Lemma~\ref{lem:sym} with $\lambda=c$ (so that $\nu_\lambda=\nu_c$), the function $F_c(y)=c\,F(\nu_cy)$ is a weak profile. By Lemma~\ref{lem:ansatz}, $\tilde u(x,t):=\tau^\beta F_c(x\tau^{-\gamma})$ is a weak solution in $\R^N\times(-\infty,0)$. Since $u_{\beta,c,s}(x,t)=\tilde u(x,t-s)$, Lemma~\ref{lem:sym} shows that $u_{\beta,c,s}$ is a weak solution in $\R^N\times(-\infty,s)$.

\emph{Value on the axis.} $u_{\beta,c,s}(0,-T)=c\,(T+s)^\beta f(0)=c\,(T+s)^\beta$.

\emph{The criterion.} Let $(x,t)\in\Pa_h$, that is, $\tau>h$ and $|x|<R(\tau)$, and let $z:=\nu_c|x|(\tau+s)^{-\gamma}$. Then
\[
c^{\frac{p-2}{p}}z=|x|\,(\tau+s)^{-\gamma}<R(\tau)(\tau+s)^{-\gamma}=\Phi_s(\tau).
\]
If $\sup_{\tau>h}\Phi_s\le y_c$, then $z<y_1$, so $f(z)>0$ and $u_{\beta,c,s}(x,t)>0$. Conversely, suppose $\Phi_s(\tau)>y_c$ for some $\tau>h$; then $y_c<\infty$. As $|x|$ ranges over $[0,R(\tau))$, the quantity $|x|(\tau+s)^{-\gamma}$ takes all the values in $[0,\Phi_s(\tau))\ni y_c$. At a point where it equals $y_c$ we have $z=y_1$ and $u_{\beta,c,s}=c(\tau+s)^\beta f(y_1)=0$, with $(x,t)\in\Pa_h$.

\emph{The cases.}
\begin{itemize}
\item If $\beta\le0$, then $y_1=\infty$ (Proposition~\ref{prop:neg}), and the condition is trivial.
\item If $0<\beta<1$, then $\gamma<\frac{p-1}{p}$ by \eqref{gamma-id}, and $\Phi_s(\tau)\sim A^{-\frac{p-1}{p}}\tau^{\frac{p-1}{p}-\gamma}\to\infty$ as $\tau\to\infty$, whereas $y_c<\infty$ by Proposition~\ref{prop:pos}. The condition fails for every $c,s,h$.
\item If $\beta=1$, then $\gamma=\frac{p-1}{p}$ and $\Phi_s(\tau)=A^{-\frac{p-1}{p}}\big(\tfrac{\tau}{\tau+s}\big)^{\frac{p-1}{p}}$, whose supremum over $\tau>h$ is $A^{-\frac{p-1}{p}}$. The condition reads $A^{-\frac{p-1}{p}}\le c^{\frac{p-2}{p}}y_1$, that is, $c\ge\big(A^{-\frac{p-1}{p}}/y_1\big)^{\frac{p}{p-2}}$. By Proposition~\ref{prop:beta1}, $y_1=k^{-\frac{p-1}{p}}$, so this is $c\ge(k/A)^{\frac{p-1}{p-2}}$, and $k^{\frac{p-1}{p-2}}=\big(N(\tfrac{p}{p-1})^{p-1}\big)^{-\frac1{p-2}}$. The formula of Remark~\ref{rem:beta1} follows from $c\,\nu_c^{\frac{p}{p-1}}=c^{\frac1{p-1}}$.
\item If $\beta>1$, then $\gamma>\frac{p-1}{p}$. For $s=0$, $\Phi_0(\tau)=A^{-\frac{p-1}{p}}\tau^{-(\gamma-\frac{p-1}{p})}$ is decreasing, so $\sup_{\tau>h}\Phi_0=A^{-\frac{p-1}{p}}h^{-(\gamma-\frac{p-1}{p})}$, which is $+\infty$ if $h=0$. The condition $\sup\Phi_0\le y_c$ reads $h\ge\big(A^{-\frac{p-1}{p}}/y_c\big)^{1/(\gamma-\frac{p-1}{p})}$. Since $y_c=c^{\frac{p-2}{p}}y_1$ and, by \eqref{gamma-id}, $\frac{p-2}{p}\cdot\frac1{\gamma-\frac{p-1}{p}}=\frac1{\beta-1}$, the right-hand side is $h_0(1)\,c^{-\frac1{\beta-1}}=h_0(c)$. For $s\ge0$ we use $\tau\le\tau+s$:
\[
\Phi_s(\tau)\ \le\ A^{-\frac{p-1}{p}}(\tau+s)^{-(\gamma-\frac{p-1}{p})}\ \le\ A^{-\frac{p-1}{p}}(h+s)^{-(\gamma-\frac{p-1}{p})}\ \le\ y_c
\]
as soon as $h+s\ge h_0(c)$. The statement of Remark~\ref{rem:beta1} follows from $h_0(c)\le h+s\iff c\ge\big(h_0(1)/(h+s)\big)^{\beta-1}$.\qed
\end{itemize}

\section{A Theorem of Phragm\'en-Lindel\"of-Type in time}
\label{section5}
Section~\ref{section4} shows that positivity in $\Pa_h$ alone does not bound the growth in time from above. Here we assume that $u$ is \emph{positive} in $\Pa_h$ and we study, for a parameter $\mu\in(0,\frac14]$ fixed below, the infimum $m_\mu(T)$ of $u(\cdot,-T)$ over a ball which is a fraction of the slice of $\Pa_h$ at time $-T$. The result is a dichotomy: either $m_\mu(T)$ is of order $T$ along a sequence $T_n\to\infty$, or $m_\mu(T)$ grows at most like $T^{\zeta_\mu}$ with an explicit exponent $\zeta_\mu<1$. The mechanism first deploys a \emph{persistence lemma}, proved by a DeGiorgi iteration. Neither the technique nor the result are intended to be new. Readers familiar with Lemma 4.1 from \cite[Chapter~3, Section~4]{DGV} will see the similarity between the mentioned lemma and Proposition \ref{prop:persist} below. We anyways include a proof of it, even if the peculiar geometry of the paraboloid is not needed. After that proof, the section purports to use this estimate as a tool for further results.

Recall the definition $R(\tau)=(\tau/A)^{\frac{p-1}{p}}$ as in \eqref{Time slice}, to mean the radius of the slice of $\Pa_h$ at time $-\tau$. As a consequence, 
\begin{equation}\label{R mu}
R(\mu\tau)=\mu^{\frac{p-1}{p}}R(\tau).
\end{equation}
Let also $\tilde\alpha:=p/N$ and $q:=p(N+2)/N=p+2\tilde\alpha$. As stated in the preliminaries, we consider an operator $\mathbf A$ which is odd in the gradient variable, and independent on $u$. As a consequence, for every $M \in \R$, if $u $ is a weak (super)solution of \eqref{p-Laplace}-\eqref{monotonicity}, then $v \coloneqq M-u$ is a (sub)solution as well.

\subsection{The persistence lemma}
In order to get a DeGiorgi iteration started, we recall the Caccioppoli energy estimate, as well as two standard tools. For the sake of brevity, we state below a particular result, which is used in the rest of the section. For a complete proof in its highest possible generality, see  \cite[Chapters~2,3]{DGV}. 

\begin{itemize}
\item \emph{Energy estimate.} Let $v$ be a weak subsolution of $v_t-\operatorname{div} \mathbf A (x,t,Dv)=0$ in $B_{2\rho}\times(t_0,t_1]$, where $\mathbf A$ satisfies \eqref{structure}, and assume $v(\cdot,t_0)\le0$. Let $k\ge0$, $w:=(v-k)_+$, and let $\chi\in C^\infty_c(B_{2\rho})$ with $0\le\chi\le1$, independent of $t$. Then
\begin{equation}
\label{lem:energy}
\esssup_{t_0<t\le t_1}\int_{B_{2\rho}}w^2\chi^p\,dx+\int_{t_0}^{t_1}\!\!\int_{B_{2\rho}}|Dw|^p\chi^p\,dx\,dt\ \le\ \Gamma_1\int_{t_0}^{t_1}\!\!\int_{B_{2\rho}}w^p|D\chi|^p\,dx\,dt,
\end{equation}
with $\Gamma_1=\Gamma_1(p,C_0,C_1)$.
\item \emph{Parabolic embedding.} If $\varphi\in L^\infty(t_0,t_1;L^2(B))\cap L^p(t_0,t_1;W^{1,p}_0(B))$ for a ball $B$, then
\begin{equation}
    \label{parabolic embedding}
    \int_{t_0}^{t_1}\!\!\int_B|\varphi|^{q}\,dx\,dt\ \le\ \Gamma_2\Big(\int_{t_0}^{t_1}\!\!\int_B|D\varphi|^p\,dx\,dt\Big)\Big(\sup_{t_0<t\le t_1}\int_B\varphi^2\,dx\Big)^{\tilde\alpha},
\end{equation}
with $\Gamma_2=\Gamma_2(N,p)$.
\item \emph{Fast geometric convergence.} Let $C>0$ and $\tilde b>1$ be constants, and let $(Y_n)$ be positive numbers satisfying 
\begin{equation}
    \label{fast geometric convergence}
    Y_{n+1}\le C\,\tilde b^{\,n}Y_n^{1+\tilde\alpha}.
\end{equation}
If $Y_0\le C^{-1/\tilde\alpha}\tilde b^{-1/\tilde\alpha^2}$, then $Y_n\to0$. Here $C$ and $\tilde b$ are generic positive constants, in the sense that the statement holds for every choice of them; of course, the same values appear in the recursion and in the smallness condition on $Y_0$.
\end{itemize}

\begin{prop}\label{prop:persist}
There is $c_0\in(0,1)$, depending only on $N,p,C_0,C_1$, such that the following holds. Let $\rho$, $S$, $M>0$ and $t_0\in\R$, let $Z:=B_{2\rho}\times(t_0,t_0+S]$, and let $u$ be a continuous weak (super)solution of \eqref{p-Laplace}-\eqref{monotonicity} in a neighborhood of $\overline Z$ with
\begin{equation}\label{persist hyp}
u\ge0\ \text{ in }Z,\qquad u(\cdot,t_0)\ge M\ \text{ in }B_{2\rho},\qquad
\varepsilon:=\frac{M^{p-2}S}{\rho^p}\ \le\ c_0 .
\end{equation}
Then $u\ge M/2$ in $B_\rho\times(t_0,t_0+S]$.
\end{prop}

The goal of the proposition is to show that a lower bound $M$ on $B_{2\rho}$ persists, as $M/2$ on $B_\rho$, for all times $s\in(0,\min\{S,c_0M^{2-p}\rho^p\}]$ after $t_0$, provided $u\ge0$ in the cylinder. The positivity of $u$ is needed only to know that $v=M-u$ is bounded above by $M$.
\begin{proof}
Let $v:=M-u$. $v$ is a weak (sub)solution of the same equation; $v\le M$ in $Z$ since $u\ge0$, and $v(\cdot,t_0)\le0$ in $B_{2\rho}$. We have to show that $v\le M/2$ in $B_\rho\times(t_0,t_0+S]$.

The argument relies on a DeGiorgi iteration. So, for $n\ge0$, let
\[
k_n:=\frac M2-\frac M{2^{n+2}},\qquad r_n:=\rho\,(1+2^{-n}),\qquad \tilde r_n:=\tfrac12(r_n+r_{n+1}),
\]
so that $k_0=M/4$, $k_n\uparrow M/2$, $r_0=2\rho$, $r_n\downarrow\rho$ and $k_{n+1}-k_n=M\,2^{-n-3}$. Let
\[
A_n:=\{(x,t)\in B_{r_n}\times(t_0,t_0+S]:\ v(x,t)>k_n\},\quad  Y_{n} := \frac{|A_n|}{|Z|}
\]
and $w_n:=(v-k_n)_+$, so that $0\le w_n\le M$ and $w_n=0$ outside $\{v>k_n\}$. Fix $\chi_n\in C^\infty_c(B_{r_n})$ with $0\le\chi_n\le1$, $\chi_n=1$ on $B_{\tilde r_n}$ and $|D\chi_n|\le 2^{n+3}/\rho$, which is possible since $r_n-\tilde r_n=\rho\,2^{-n-2}$.

We now apply the energy estimate \eqref{lem:energy} with $k=k_n$, $\chi=\chi_n$, together with $w_n\le M$ and $w_n=0$ outside $\{v>k_n\}$. We then get
\begin{equation}\label{X n}
\esssup_{t}\int_{B_{2\rho}}w_n^2\chi_n^p\,dx+\iint_{Z}|Dw_n|^p\chi_n^p\,dx\,dt\ \le\ \Gamma_1\,\frac{2^{(n+3)p}M^p}{\rho^p}\,|A_n|=:X_n .
\end{equation}
Set $\varphi_n:=w_n\chi_n^{p/2}$. $\varphi_n$ is supported in $B_{r_n}$, and $\varphi_n \in L^\infty(t_0,t_0+S;L^2)\cap L^p(t_0,t_0+S;W^{1,p}_0(B_{2\rho}))$. Moreover, $\esssup_t\int\varphi_n^2\le X_n$ by \eqref{X n}; and since $0\le\chi_n\le1$ and $p\ge2$,
via Young and standard inequalities we have
\[
\begin{gathered}
|D\varphi_n|^p\le2^{p-1}\Big(\chi_n^{p}|Dw_n|^p+\big(\tfrac p2\big)^pw_n^p|D\chi_n|^p\Big),\\
\iint_Z w_n^p|D\chi_n|^p\le\frac{2^{(n+3)p}M^p}{\rho^p}|A_n|=\frac{X_n}{\Gamma_1},
\end{gathered}
\]
so that $\iint_Z|D\varphi_n|^p\le\Gamma_3X_n$ with $\Gamma_3=\Gamma_3(p,C_0,C_1)$.

On $A_{n+1}$ we have $|x|<r_{n+1}<\tilde r_n$, hence $\chi_n=1$, and $w_n>k_{n+1}-k_n=M\,2^{-n-3}$. Therefore, by the parabolic embedding \eqref{parabolic embedding} for $B_{2\rho}\times(t_0,t_0+S]$,
\[
\Big(\frac M{2^{n+3}}\Big)^{q}|A_{n+1}|\ \le\ \iint_{A_{n+1}}\varphi_n^{\,q}\ \le\ \iint_Z\varphi_n^{\,q}\ \le\ \Gamma_2\Gamma_3\,X_n^{\,1+\tilde\alpha}.
\]
We insert \eqref{X n} and divide by $|Z|=|B_1|(2\rho)^NS$. As a result, 
\begin{equation}\label{eq:rec}
Y_{n+1}\ \le\ C_*\,\tilde b^{\,n}\ \varepsilon^{\tilde\alpha}\,Y_n^{1+\tilde\alpha},\qquad\varepsilon=\frac{M^{p-2}S}{\rho^p},
\end{equation}
with $\tilde b:=2^{q+p(1+\tilde\alpha)}>1$ and $C_*=C_*(N,p,C_0,C_1)>0$. The coefficient depends on $M,\rho,S$ only through $\varepsilon$.

We now apply \eqref{fast geometric convergence} with $C$ replaced by $C_*\varepsilon^{\tilde\alpha}$. So, $Y_n\to0$ provided
\[
Y_0\ \le\ \big(C_*\varepsilon^{\tilde\alpha}\big)^{-1/\tilde\alpha}\tilde b^{-1/\tilde\alpha^2}=C_*^{-1/\tilde\alpha}\,\tilde b^{-1/\tilde\alpha^2}\,\varepsilon^{-1}.
\]
Since $Y_0\le1$ by its definition, this holds if $\varepsilon\le c_0:=\min\{1/2,\ C_*^{-1/\tilde\alpha}\tilde b^{-1/\tilde\alpha^2}\}$. Then $v\le M/2$ a.e., and everywhere by continuity.
\end{proof}

\subsection{The scale step}\label{sec:scale}
We now aim to start from the previous estimate in order to obtain a chain of inequalities connecting different time levels $-T < -t < -h$. The argument relies on the following quantities, $m_\mu(T)$ and $\tilde Q_\mu(T)$, to be defined with a given parameter $\mu$. 

Fix $\mu\in(0,\frac14]$ and define $\delta_\mu\in(0,1)$ by $\delta_\mu^{\frac{p}{p-1}}=\mu/2$. We stress that $\delta_\mu$ is not a universal constant like the $\delta$ of the expansion of positivity: it is just a number, depending only on $\mu$ and $p$, which fixes the size of the ball on which $m_\mu$ is computed. Let
\begin{equation*}
m_\mu(T):=\inf_{|x|<\delta_\mu R(T)}u(x,-T),\qquad \tilde Q_\mu(T):=\frac{m_\mu(T)}{T}\qquad(T>h).
\end{equation*}
Since $\delta_\mu<1$, the closed ball $\overline{B_{\delta_\mu R(T)}}$ is a compact subset of the slice of $\Pa_h$; so $m_\mu$ is positive because $u$ is. With $c_0$ from Proposition~\ref{prop:persist}, set
\begin{equation}\label{eq:qmu}
\kappa_0:=\big(c_0\,2^{-(2p-1)}\big)^{\frac1{p-2}},\qquad
q_\mu:=\kappa_0\,A^{-\frac{p-1}{p-2}}\,\mu^{\frac{p-1}{p-2}} .
\end{equation}

\begin{prop}\label{prop:step}
Let $\mu\in(0,\frac14]$ and $T>h/\mu$. If $\tilde Q_\mu(T)\le q_\mu$, then $m_\mu(\mu T)\ge\frac12\,m_\mu(T)$.
\end{prop}
\begin{proof}
Put $M:=m_\mu(T)$, $\rho:=\frac12\delta_\mu R(T)$, $t_0:=-T$ and $S:=(1-\mu)T$, so that $B_{2\rho}=B_{\delta_\mu R(T)}$ and $t_0+S=-\mu T$.

\emph{The cylinder $Z:=B_{2\rho}\times(-T,-\mu T] $ lies in $\Pa_h$.} For $\tau\in[\mu T,T]$, \eqref{R mu} and $\delta_\mu^{\frac p{p-1}}=\mu/2<\mu$ give $\delta_\mu R(T)<\mu^{\frac{p-1}{p}}R(T)=R(\mu T)\le R(\tau)$. As $\mu T>h$, we get $\overline{B_{2\rho}}\times[-T,-\mu T]\subset\Pa_h$. Hence $u>0$ in the cylinder $Z$, and $u(\cdot,-T)\ge M$ in $B_{2\rho}$ by the definition of $m_\mu(T)$.

\emph{The condition $\tilde Q_\mu(T)\le q_\mu$.} We have $\rho^p=2^{-p}\delta_\mu^pR(T)^p=2^{-p}(\mu/2)^{p-1}(T/A)^{p-1}$, and $M=\tilde Q_\mu(T)\,T$. Since $S\le T$, we want to check whether $\varepsilon=M^{p-2}S/\rho^p\le c_0$. It suffices that  
\[
1 \leq \frac{c_0\,M^{2-p}\rho^p}{T}=c_0\,2^{-(2p-1)}\mu^{p-1}A^{-(p-1)}\,\tilde Q_\mu(T)^{-(p-2)}\ 
\quad\Longleftrightarrow\quad \tilde Q_\mu(T)\le q_\mu .
\]
Now, Proposition~\ref{prop:persist} gives $u\ge M/2$ in $B_\rho\times(-T,-\mu T]$. At the time $-\mu T$, $\delta_\mu R(\mu T)=\mu^{\frac{p-1}{p}}\delta_\mu R(T)\le\frac12\delta_\mu R(T)=\rho$, since $\mu\le\frac14$ and $\frac{p-1}p>\frac12$. Hence we conclude $m_\mu(\mu T)\ge M/2$.
\end{proof}

\begin{cor}\label{cor:iter}
Let $\mu\in(0,\frac14]$, $\tau_0>h$ and $T_i:=\mu^{-i}\tau_0$. If $\tilde Q_\mu(T_i)\le q_\mu$ for every $i=1,\dots,j$, then
\[
m_\mu(T_j)\ \le\ \Big(\frac{T_j}{\tau_0}\Big)^{\zeta_\mu}m_\mu(\tau_0),\qquad \zeta_\mu:=\frac{\ln2}{\ln(1/\mu)}\ \le\ \frac12 .
\]
\end{cor}
\begin{proof}
Since $\mu T_i=T_{i-1}>h$, Proposition~\ref{prop:step} at time $T_i$ gives $m_\mu(T_{i-1})\ge\frac12m_\mu(T_i)$, that is $m_\mu(T_i)\le2\,m_\mu(T_{i-1})$. Hence $m_\mu(T_j)\le2^jm_\mu(\tau_0)$, and $2^j=\big(\mu^{-j}\big)^{\zeta_\mu}=(T_j/\tau_0)^{\zeta_\mu}$ by the definition of $\zeta_\mu$.
\end{proof}

\subsection{The dichotomy}
We now state the following theorem. The two alternatives below are complementary. 
\begin{thm}\label{thm:dich}
Let $u$ be a local weak solution to \eqref{p-Laplace}-\eqref{monotonicity} in $\R^{N} \times (-\infty, 0]$, and suppose that $u>\Theta$ in $\mathcal{P}_h$ for some $A>0$, $\Theta\in\R$, $h\ge0$. Fix $\mu\in(0,\frac14]$ and define $m_\mu$ with $u-\Theta$ in place of $u$. Then, one of the following holds:
\begin{enumerate}
\item[(i)] $\limsup_{T\to\infty}\tilde Q_\mu(T)\ge q_\mu$. In particular, for every $\bar q<q_\mu$ there exist $T_n\to\infty$ with $m_\mu(T_n)\ge \bar q\,T_n$.
\item[(ii)] $\limsup_{T\to\infty}\tilde Q_\mu(T)<q_\mu$. Then $m_\mu(T)=O\big(T^{\zeta_\mu}\big)$ as $T\to\infty$, with $\zeta_\mu=\ln2/\ln(1/\mu)\le\frac12$.
\end{enumerate} 
In particular, if $\tilde Q_\mu(T)\to0$ for every $\mu\in(0,\frac14]$, then for every $\sigma>0$ there is $\mu=\mu(\sigma)$ such that $m_\mu(T)=O(T^{\sigma})$; note that the ball $B_{\delta_\mu R(T)}$ shrinks as $\mu\to0$.
\end{thm}
\begin{proof}
Without loss of generality, let $\Theta=0$. Assume first alternative (i). For $\bar q<q_\mu$ there is a sequence of $T_n\to\infty$ with $\tilde Q_\mu(T_n)\ge \bar q$, which is the claim.

Assume then alternative (ii). Let $T_\star>h/\mu$ be such that $\tilde Q_\mu(T)\le q_\mu$ for every $T\ge T_\star$. Any $T\ge T_\star$ can be written as $T=\mu^{-j}\tau_0$, for some $j\ge0$ and $\tau_0\in[T_\star,T_\star/\mu)$. For every $i\le j$, the time instant $T_i=\mu^{-i}\tau_0 \ge T_\star$, so Corollary~\ref{cor:iter} applies and
\[
m_\mu(T)\ \le\ \Big(\frac{T}{\tau_0}\Big)^{\zeta_\mu}m_\mu(\tau_0)\ \le\ T_\star^{-\zeta_\mu}\Big(\max_{[T_\star,T_\star/\mu]}m_\mu\Big)\,T^{\zeta_\mu},
\]
and the maximum is finite since $m_\mu$ is continuous (or directly, $m_\mu(T)\le u(0,-T)$). In other words, $m_\mu(T)\le C_{T_\star}T^{\zeta_\mu}$ for every $T\ge T_\star$, where
\[
C_{T_\star}:=T_\star^{-\zeta_\mu}\max_{[T_\star,T_\star/\mu]}m_\mu
\]
depends on $u$, $\mu$ and $T_\star$. If $\tilde Q_\mu(T)\to0$ for every $\mu$, then (ii) holds for every $\mu\in(0,\frac14]$, and $\zeta_\mu\to0$ as $\mu\to0$.
\end{proof}

\subsection{From the infimum to the supremum}
In this subsection we show that, in alternative (ii) of Theorem~\ref{thm:dich}, a Harnack-type argument turns the bound on the infimum into a bound on the supremum over the same ball. Let $C_H,\Gamma_H$ be the constants of the degenerate Harnack inequality (Theorem~\ref{thm:harnack}), which depend only on the data. Assume $u>0$ in $\Pa_h$ (otherwise replace $u$ by $u-\Theta$), and set
\begin{equation}\label{eq:Mmu}
M_\mu(T):=\sup_{|x|<\delta_\mu R(T)}u(x,-T),
\end{equation}
which is finite since $\overline{B_{\delta_\mu R(T)}}\subset\Pa_h$ is compact and $u$ is continuous. Let
\begin{equation}\label{eq:mu0c1}
\bar\mu:=\big(1+8\sqrt N\big)^{-\frac p{p-1}},\qquad 
c_1:=C_H\,2^{\frac1{p-2}}\,(8\delta_\mu)^{\frac p{p-2}}A^{-\frac{p-1}{p-2}} .
\end{equation}
Note that $\bar\mu<\frac14$, and that, since $\delta_\mu^{\frac{p}{p-2}}=(\mu/2)^{\frac{p-1}{p-2}}$ (we recall that $\delta_\mu$ is the number of Section~\ref{sec:scale}, not the constant $\delta$ of the expansion of positivity), the constant $c_1$ has the same dependence on $\mu$ and $A$ as $q_\mu$, namely $c_1\simeq(\mu/A)^{\frac{p-1}{p-2}}$ with a factor depending only on the data.

\begin{lem}[From $\inf$ to $\sup$]\label{lem:harnack}
Let $u$ be as in Theorem~\ref{thm:dich} with $\Theta=0$, let $\mu\le\bar\mu$ and $T>2h$. Then
\[
M_\mu(T)\ \le\ \max\Big\{c_1T,\ \Gamma_H\sup_{T/2\le S\le T}m_\mu(S)\Big\}.
\]
\end{lem}
\begin{proof}
Fix $x$ with $|x|<\delta_\mu R(T)$, let $u_0:=u(x,-T)>0$ and assume $u_0\ge c_1T$ (otherwise there is nothing to prove). Set
\[
\rho:=|x|+\delta_\mu R(T)<2\delta_\mu R(T),\qquad \theta:=(C_H/u_0)^{p-2}.
\]
\emph{Step 1: the time lag.} Since $(4\rho)^p<(8\delta_\mu)^pR(T)^p=(8\delta_\mu)^p(T/A)^{p-1}$,
\[
\theta(4\rho)^p< C_H^{p-2}u_0^{2-p}(8\delta_\mu)^pA^{-(p-1)}T^{p-1}\le\frac T2,
\]
where the last inequality is equivalent to $u_0\ge c_1T$. Hence $T'':=T-\theta\rho^p\in[T/2,T]$, and $T''>h$ because $T>2h$.

\emph{Step 2: the cylinder lies in $\Pa_h$.} Consider the intrinsic cylinders $(x,-T)+Q_{4\rho}^{\pm}(\theta)$. Their time levels $-\tau$ satisfy $\tau\ge T-\theta(4\rho)^p\ge T/2>h$. Moreover $K_{4\rho}(x)\subset B_{4\sqrt N\rho}(x)$, so every point $y$ of the cylinders satisfies
\[
|y|<|x|+4\sqrt N\rho<\delta_\mu\big(1+8\sqrt N\big)R(T)\le 2^{-\frac{p-1}{p}}R(T)=R(T/2)\le R(\tau),
\]
where the inequality in the middle is equivalent to $\mu\le\bar\mu$, because $\delta_\mu^{\frac p{p-1}}=\mu/2$. Hence the cylinders are contained in $\Pa_h$.

\emph{Step 3: Harnack.} By Step~2, $u$ is a continuous positive solution in the cylinder $K_{4\rho}(x)\times(-T-\theta(4\rho)^p,-T+\theta(4\rho)^p]\subset\Pa_h$, so the second inequality of the Harnack inequality (after a time translation) gives
\[
u_0\le\Gamma_H\inf_{K_\rho(x)}u(\cdot,-T'').
\]
Since $T''\le T$, for $|y|<\delta_\mu R(T'')$ we have $|y-x|<\delta_\mu R(T)+|x|=\rho$, so $B_{\delta_\mu R(T'')}\subset B_\rho(x)\subset K_\rho(x)$. Therefore $\inf_{K_\rho(x)}u(\cdot,-T'')\le m_\mu(T'')$, and $u_0\le\Gamma_H\,m_\mu(T'')$ with $T''\in[T/2,T]$.
\end{proof}

\begin{thm}[Linear growth of $\sup$ in alternative (ii)]\label{thm:sup}
Let $u$ be as in Theorem~\ref{thm:dich}, $\mu\le\bar\mu$, and assume that alternative (ii) of Theorem~\ref{thm:dich} holds. Then there is $T_1>0$ such that
\[
\sup_{|x|<\delta_\mu R(T)}\big(u(x,-T)-\Theta\big)\le c_1\,T\qquad\text{for all }T\ge T_1 .
\]
\end{thm}
\begin{proof}
Without loss of generality $\Theta=0$. Let $T_\star$ and $C_{T_\star}$ be as in the proof of Theorem~\ref{thm:dich}, so that $m_\mu(T)\le C_{T_\star}T^{\zeta_\mu}$ for $T\ge T_\star$. For $T\ge2T_\star$ we have $T>2h$ and $S\ge T_\star$ for every $S\in[T/2,T]$, so Lemma~\ref{lem:harnack} gives $M_\mu(T)\le\max\{c_1T,\Gamma_HC_{T_\star}T^{\zeta_\mu}\}$. Since $\zeta_\mu<1$, we have $\Gamma_HC_{T_\star}T^{\zeta_\mu}\le c_1T$ as soon as $T\ge(\Gamma_HC_{T_\star}/c_1)^{1/(1-\zeta_\mu)}$. Take $T_1$ as the maximum of this number and $2T_\star$.
\end{proof}

\begin{rem}
The linear threshold in Lemma~\ref{lem:harnack} comes from intrinsic scaling: $\theta\rho^p\lesssim T$ if and only if $u\gtrsim T$, when $\rho\simeq R(T)$. So this argument gives $\sup\lesssim T$ and not $\sup=O(T^{\sigma})$ with $\sigma<1$, even when $\tilde Q_\mu\to0$ for every $\mu$. It controls the supremum over the same ball $B_{\delta_\mu R(T)}$ as $m_\mu$: a larger ball would force a larger $\rho$ and the Harnack cylinders would leave $\Pa_h$. We do not know whether $\sup=o(T)$ can occur in alternative (ii).
\end{rem}

\begin{rem}\label{rem:picture}
Fix $\mu\le\frac14$. Concerning the threshold $q_\mu$, we may say that 
\begin{itemize}
\item[$\bullet$] $q_\mu=\kappa_0A^{-\frac{p-1}{p-2}}\mu^{\frac{p-1}{p-2}}$ increases with $\mu$ and depends on the equation only through $c_0$;
\item[$\bullet$] below the threshold, $\zeta_\mu=\ln2/\ln(1/\mu)$ tends to $0$ as $\mu\to0$, and equals $\frac12$ for $\mu=\frac14$;
\item[$\bullet$] above the threshold, we only know that $m_\mu$ is comparable to $T$ along a sequence.
\end{itemize}
\end{rem}

\begin{rem}\label{rem:consistency}
As a consistency check with the self-similar solutions, we compare Theorem~\ref{thm:dich} with the solutions $u_{\beta,c,0}$ of Theorem~\ref{thm:positivity}, for $\mathbf A(x,t,\xi)=|\xi|^{p-2}\xi$.  Here $\mu\in(0,\frac14]$ is fixed, $\delta_\mu$ is the number of Section~\ref{sec:scale} (recall that it is not the constant $\delta$ of the expansion of positivity), and $c_0=c_0(N,p)$ is the constant of Proposition~\ref{prop:persist} for $C_0=C_1=1$.
 
\emph{The cases $\beta\le0$ and $\beta>1$.} If $\beta\le0$, then $f$ is nondecreasing (Theorem~\ref{thm:profile}), so the infimum over the ball is attained at $x=0$ and $m_\mu(T)=u(0,-T)=c\,T^{\beta}=O(1)$: alternative (ii) of Theorem~\ref{thm:dich} holds. If $\beta>1$, then $\gamma>\frac{p-1}{p}$ and the argument $z=\nu_c\,\delta_\mu R(T)\,T^{-\gamma}$ of $f$ tends to $0$ as $T\to\infty$. Hence $f(z)\to1$ and $m_\mu(T)\sim c\,T^{\beta}$, so $\tilde Q_\mu(T)\to\infty$ and alternative (i) of Theorem~\ref{thm:dich} holds. For $0<\beta<1$ there is no solution in the family which is positive in $\Pa_h$ (Theorem~\ref{thm:positivity}), so there is nothing to compare.
 
\emph{The case $\beta=1$: the values of $m_\mu$ and $\tilde Q_\mu$.} Let $c\ge a(A)$, so that $u=u_{1,c,0}$ is positive in $\Pa_h$, and recall $\gamma=\frac{p-1}{p}$ and $f(z)=1-k\,z^{\frac p{p-1}}$. We proceed in four steps.
\begin{enumerate}
\item[\emph{Step 1.}] \emph{$m_\mu$ is linear in $T$.} Since $R(T)=A^{-\frac{p-1}{p}}T^{\gamma}$, for $|x|<\delta_\mu R(T)$ the argument $z=\nu_c|x|T^{-\gamma}$ of $f$ satisfies
\[
z\ <\ z_*:=\nu_c\,\delta_\mu\,A^{-\frac{p-1}{p}}=c^{-\frac{p-2}{p}}\,\delta_\mu\,A^{-\frac{p-1}{p}},
\]
and $z_*$ does not depend on $T$. As $f$ is decreasing, the infimum of $u(\cdot,-T)=c\,T\,f(z)$ over the ball is the limit at the boundary: $m_\mu(T)=c\,T\,f(z_*)$. Using $\delta_\mu^{\frac p{p-1}}=\mu/2$,
\[
z_*^{\frac p{p-1}}=\frac\mu2\,\frac1A\,c^{-\frac{p-2}{p-1}},
\qquad\text{so}\qquad
\tilde Q_\mu(T)=c\,f(z_*)=c\Big(1-\frac\mu2\,\vartheta_c\Big)=:\tilde Q_c,
\]
where $\vartheta_c:=\frac kA\,c^{-\frac{p-2}{p-1}}$.
The number $\vartheta_c$ is decreasing in $c$ and, because $a(A)=(k/A)^{\frac{p-1}{p-2}}$, we have $\vartheta_c\le1\iff c\ge a(A)$, with $\vartheta_{a(A)}=1$. In particular $\tilde Q_c\ge c(1-\frac\mu2)>0$. Note that $\tilde Q_c=c-\frac\mu2\frac kA\,c^{\frac1{p-1}}$.
\item[\emph{Step 2.}] \emph{Alternative (i) of Theorem~\ref{thm:dich} holds.} Since $m_\mu(T)=\tilde Q_c\,T$ with $\tilde Q_c>0$, $m_\mu$ is not $O(T^{\zeta_\mu})$, because $\zeta_\mu\le\frac12<1$. So alternative (ii) fails and alternative (i) holds, that is $\limsup Q=\tilde Q_c\ge q_\mu$. This is true for every $c\ge a(A)$.
\item[\emph{Step 3.}] \emph{The smallest $\tilde Q_c$.} We have
\[
\frac{d}{dc}\tilde Q_c=1-\frac{\mu}{2(p-1)}\,\vartheta_c\ \ge\ 1-\frac{\mu}{2(p-1)}>0\qquad(c\ge a(A)),
\]
so $c\mapsto \tilde Q_c$ is increasing on $[a(A),\infty)$ and its minimum is $\tilde Q_{a(A)}=(1-\frac\mu2)\,a(A)$. By Step~2,
\[
q_\mu\ \le\ \Big(1-\frac\mu2\Big)\,a(A)\qquad\text{for every }\mu\in(0,\tfrac14].
\]
\item[\emph{Step 4.}] \emph{What it says on $c_0$.} By \eqref{eq:qmu} and $a(A)=k^{\frac{p-1}{p-2}}A^{-\frac{p-1}{p-2}}$, the power of $A$ cancels. Raising to the power $p-2$, the inequality of Step~3 becomes
\[
c_0\,2^{-(2p-1)}\,\mu^{p-1}\ \le\ \Big(1-\frac\mu2\Big)^{p-2}k^{p-1},\qquad k^{p-1}=\frac1N\Big(\frac{p-1}{p}\Big)^{p-1}.
\]
The left-hand side increases with $\mu$ and the right-hand side decreases, so the strongest case is $\mu=\frac14$, which gives $c_0\le2^{4p-3}\big(\tfrac78\big)^{p-2}k^{p-1}$. This is not restrictive, since already $c_0\le\frac12$ and the right-hand side exceeds $1$ unless $N$ is very large; it is only a consistency check.
\end{enumerate}
\end{rem}

\section{The singular supercritical case: the result in space}\label{section6}
From now on, up to the end of Section~\ref{section9}, we leave the degenerate range and consider the \emph{singular, supercritical} range
\[
p_*:=\frac{2N}{N+1}<p<2,\qquad p':=\frac p{p-1}>2 .
\]
The assumption $p>p_*$ is exactly the one under which the singular Harnack inequality of Section~\ref{section2} is available; in particular, by Remark~\ref{Rmk-Hoelder} the solutions we deal with are locally H\"older continuous, and hence continuous. Sections~\ref{section6}--\ref{section9} are the singular counterparts of Sections~\ref{section3}--\ref{section5}, in the same order: the result in space (this section, and Section~\ref{section7} for the case in which the iteration cannot be started), the self-similar profiles (Section~\ref{section8}) and the result in time (Section~\ref{section9}). The equation is \eqref{p-Laplace}-\eqref{monotonicity} with the structure conditions of Section~\ref{sec:equations}; the Comparison Principle and Remark~\ref{rem:shift} hold as before, and the expansion of positivity we use is Proposition~\ref{exp pos sing}. The paraboloid $\Pa_h$ is the set \eqref{Pa}, whose exponent $\frac p{p-1}$ is now denoted by $p'$, and $R(\tau)$ is as in \eqref{Time slice}. The symbols $\xi_\eta,\omega,C_\eta,\Lambda$ introduced below are redefined here for the singular case, and are not those of the proof of Theorem~\ref{P-L Thm in space}.

Let $\eta,\delta$ be the constants of Proposition~\ref{exp pos sing} with $\alpha=\frac12$ and put
\begin{gather*}
\xi_\eta:=\frac\eta{1-\eta},\quad \omega:=\delta\,16^{-p},\quad C_\eta:=\Big(1-(1+\xi_\eta)^{-\frac{p-1}p}\Big)^{-1},\\
\Lambda:=1+\sqrt N\,C_\eta, \quad \Xi_*:=\Big(\frac{A\Lambda^{p'}}{\omega}\Big)^{\frac1{2-p}},\quad L_*:=\Xi_*\Lambda^{-p'}=\Big(\frac{A\Lambda^{p}}{\omega}\Big)^{\frac1{2-p}} .
\end{gather*}
Note that $\omega\,\Xi_*^{2-p}=A\Lambda^{p'}$.

\begin{thm}\label{s:thm:space}
Let $u$ be a local weak solution of \eqref{p-Laplace}-\eqref{monotonicity} in $\R^N\times(-\infty,0]$ and suppose $u\le\Theta$ in $\Pa_h$. Then one of the following holds.
\begin{enumerate}
\item[(i)] For every $(x,t)$ with $t<-h$ and $u(x,t)>\Theta$ one has $u(x,t)-\Theta<\Xi_*|x|^{p'}$. In particular $u(0,t)\le\Theta$ for $t<-h$.
\item[(ii)] There is $(x_0,t_0)$, $t_0<-h$, with $u(x_0,t_0)>\Theta$ and $u(x_0,t_0)-\Theta\ge\Xi_*|x_0|^{p'}$, and
\[
\limsup_{\rho\to\infty}\ \sup_{|x|=\rho,\ t<-h}\frac{u(x,t)-\Theta}{\rho^{p'}}\ \ge\ L_*>0 .
\]
\end{enumerate}
\end{thm}

\begin{proof}
By Remark~\ref{rem:shift} we may take $\Theta=0$. Suppose (i) fails: there is $(x_0,t_0)$, $t_0<-h$, with $u_0:=u(x_0,t_0)>0$ and $u_0\ge\Xi_*|x_0|^{p'}$.

\smallskip\noindent\emph{Construction.} We build $(x_k,t_k)$, $u_k:=u(x_k,t_k)>0$, $t_{k+1}\le t_k$, and set
\[
\rho_k:=\Big(\frac{u_k}{\Xi_*}\Big)^{\frac1{p'}},\quad s_k:=t_k-\omega\rho_k^pu_k^{2-p},\quad Z^k:=K_{\rho_k}(x_k)\times(s_k,t_k],\quad M_k:=\max_{\overline{Z^k}}u .
\]
Since $p+p'(2-p)=p'$ we have $\rho_k^pu_k^{2-p}=\Xi_*^{2-p}\rho_k^{p'}$, hence
\begin{equation}\label{s:eq:sk}
s_k=t_k-\omega\Xi_*^{2-p}\rho_k^{p'}=t_k-A\Lambda^{p'}\rho_k^{p'} .
\end{equation}
By the choice of $(x_0,t_0)$ we have $\rho_0\ge|x_0|$; this is the only place where the hypothesis $u_0\ge\Xi_*|x_0|^{p'}$ enters. We shall prove that one can choose $(x_{k+1},t_{k+1})\in\overline{Z^k}$ with
\begin{equation}\label{s:eq:growth2}
u_{k+1}\ge(1+\xi_\eta)u_k,\qquad |x_{k+1}-x_k|_\infty=\rho_k,\qquad t_{k+1}\le t_k .
\end{equation}

\smallskip\noindent\emph{Step 0: geometry.} Assume \eqref{s:eq:growth2} up to $k$. Then $u_j\le(1+\xi_\eta)^{-(k-j)}u_k$, hence $\rho_j\le(1+\xi_\eta)^{-(k-j)\frac{p-1}p}\rho_k$, and $|x_0|\le\rho_0\le\rho_k$. For $x\in K_{\rho_k}(x_k)$,
\begin{equation}\label{s:eq:xb}
|x|\le|x_0|+\sqrt N\sum_{j\le k}\rho_j\le\big(1+\sqrt NC_\eta\big)\rho_k=\Lambda\rho_k .
\end{equation}

\smallskip\noindent\emph{Step 1: $\overline{K_{\rho_k}(x_k)}\times\{s_k\}\subset\Pa_h$.} Since $t_k<0$, \eqref{s:eq:sk} and \eqref{s:eq:xb} give $A|x|^{p'}\le A\Lambda^{p'}\rho_k^{p'}<-s_k$ for $x\in\overline{K_{\rho_k}(x_k)}$; and $s_k<t_k\le t_0<-h$. So $u\le0$ there.

\smallskip\noindent\emph{Step 2: choice of the next vertex.} By the Comparison Principle $M_k$ is attained on $\partial_pZ^k$. On the bottom $u\le0<u_k\le M_k$, so $M_k=u(x_{k+1},t_{k+1})$ with $|x_{k+1}-x_k|_\infty=\rho_k$ and $s_k<t_{k+1}\le t_k$.

\smallskip\noindent\emph{Step 3: $M_k\ge(1+\xi_\eta)u_k$.} Define $r_k$ by
\[
\delta M_k^{2-p}r_k^p=\omega\rho_k^pu_k^{2-p},\qquad\text{i.e.}\qquad r_k^p=16^{-p}\rho_k^p\Big(\frac{u_k}{M_k}\Big)^{2-p}\le\Big(\frac{\rho_k}{16}\Big)^p,
\]
because $M_k\ge u_k$. Hence $K_{16r_k}(x_k)\times(s_k,t_k]\subset Z^k$, the function $v:=M_k-u$ is a nonnegative solution in $Z^k$ (Remark~\ref{rem:shift}), and $v(\cdot,s_k)\ge M_k$ in $K_{r_k}(x_k)$ by Step 1. Proposition~\ref{exp pos sing}, applied in $E_T=Z^k$ with $\alpha=\frac12$ (the datum holds with $\alpha=1$, hence with $\frac12$), $M=M_k$, $\rho=r_k$, gives $v(\cdot,t)\ge\eta M_k$ at $x_k$ for
\[
s_k+(1-\bar\eps)\,\omega\rho_k^pu_k^{2-p}\le t\le s_k+\omega\rho_k^pu_k^{2-p}=t_k .
\]
The time $t_k$ is the right endpoint of the window, so $M_k-u_k=v(x_k,t_k)\ge\eta M_k$, that is $M_k\ge\frac{u_k}{1-\eta}=(1+\xi_\eta)u_k$.

\smallskip\noindent\emph{Conclusion.} The sequence is well defined and $u_k\ge(1+\xi_\eta)^ku_0\to\infty$, so $\rho_k\to\infty$. Since $\rho_k\le|x_{k+1}-x_k|\le|x_k|+|x_{k+1}|$, we have $|x_k|\to\infty$ along a subsequence. By \eqref{s:eq:xb} at $x=x_k$,
\[
u_k=\Xi_*\rho_k^{p'}\ge\Xi_*\Lambda^{-p'}|x_k|^{p'}=L_*|x_k|^{p'},\qquad t_k\le t_0<-h,
\]
which is (ii).
\end{proof}

\begin{rem}\label{s:rem:slope}
Along the construction the quotient $\Xi:=u_k/\rho_k^{p'}$ is constant, and Step~1 is exactly the requirement $\omega\,\Xi^{2-p}\ge A\Lambda^{p'}$, that is $\Xi\ge\Xi_*$: the \emph{slope} of the starting point must be at least $\Xi_*$. Any $\Xi\ge\Xi_*$ would work, with $\rho_0=(u_0/\Xi)^{1/p'}$, but the condition $\rho_0\ge|x_0|$ reads $u_0\ge\Xi|x_0|^{p'}$ and is weakest for $\Xi=\Xi_*$. This is why the threshold $u_0\ge\Xi_*|x_0|^{p'}$ appears, and why the alternative of the theorem is stated in terms of the cone $\{u-\Theta<\Xi_*|x|^{p'}\}$.
\end{rem}

\section{The case $u_0<\Xi_*|x_0|^{p'}$}\label{section7}

Theorem~\ref{s:thm:space}(i) is the ``bad'' case, in which the iteration cannot be started. The following explicit family shows that it cannot be excluded.

\begin{prop}[Explicit solutions]\label{s:prop:family}
For $c>0$ let $\kappa_c:=N(cp')^{p-1}$ and $w_c(x,t):=c|x|^{p'}+\kappa_ct$. Then $w_c$ is an eternal solution of $u_t-\Delta_pu=0$. Moreover, with
\[
c_A:=\big(NA\,p'^{\,p-1}\big)^{\frac1{2-p}},
\]
one has $w_c\le0$ in $\Pa_h$ if and only if $c\le c_A$. For such $c$, $w_c>0$ at points with $t<-h$ outside $\Pa_h$, and
\[
\limsup_{\rho\to\infty}\sup_{|x|=\rho,\ t<-h}\frac{w_c}{\rho^{p'}}=c .
\]
\end{prop}
\begin{proof} It is a matter of direct computation to see that 
$Dw_c=cp'|x|^{p'-2}x$ and $|Dw_c|^{p-2}Dw_c=(cp')^{p-1}|x|^{(p'-1)(p-2)+p'-2}x=(cp')^{p-1}x$, since $(p'-1)(p-2)+p'-2=0$. Hence $\Delta_pw_c=N(cp')^{p-1}=\kappa_c=\partial_tw_c$.

In $\Pa_h$ we have $t<\min\{-A|x|^{p'},-h\}$. For $|x|^{p'}\ge h/A$ the supremum of $w_c(x,\cdot)$ over this range is $(c-\kappa_cA)|x|^{p'}$; for $|x|^{p'}<h/A$ it is $c|x|^{p'}-\kappa_ch<h(c/A-\kappa_c)$. Hence $w_c\le0$ in $\Pa_h$ if and only if $c\le\kappa_cA=NAp'^{\,p-1}c^{p-1}$, i.e. $c^{2-p}\le NAp'^{\,p-1}$ (here $p<2$ is used).

For $c\le c_A$ we have $c/\kappa_c\le A$, so $w_c>0$ on $\{-\frac c{\kappa_c}|x|^{p'}<t<-h\}$; this set is nonempty for $|x|$ large and lies outside $\Pa_h$. Finally $\sup_{|x|=\rho,\,t<-h}w_c/\rho^{p'}=c-\kappa_ch\rho^{-p'}\to c$.
\end{proof}

\begin{rem}\label{s:rem:consequences}
For every $c\in(0,c_A]$, $w_c$ satisfies the hypotheses of Theorem~\ref{s:thm:space} with $\Theta=0$, is positive at points of $\{t<-h\}$, and has $\limsup=c$, with $c$ arbitrarily small. Since $w_c<c|x|^{p'}$, for $c<\Xi_*$ it falls in alternative (i) of Theorem~\ref{s:thm:space}. So in the bad case there is no universal positive lower bound for the $\limsup$, and the quantitative statement of Theorem~\ref{s:thm:space}, with the threshold $\Xi_*$, is what this method gives.
\end{rem}

\section{Self-similar profiles for $p_*<p<2$}\label{section8}

Throughout this section $p_*<p<2$. We keep the notation of Section~\ref{section4}: $\phi_p(s)=|s|^{p-2}s$, $\psi=\phi_p^{-1}$, $\gamma=\gamma_\beta$ as in \eqref{gamma}, and the identities \eqref{gamma-id}, which hold for every $p$. For $p<2$ we write
\[
\psi(s)=\operatorname{sign}(s)|s|^{\varsigma},\qquad\varsigma:=\frac1{p-1}>1 .
\]
Thus $\psi\in C^1(\R)$ is locally Lipschitz (it was not Lipschitz at $0$ for $p>2$), but it is superlinear at infinity. Put
\[
\beta_c:=\frac1{2-p}>1,\qquad\varkappa:=N-(N-1)\varsigma .
\]
Then $\gamma>0\iff\beta<\beta_c$, and \eqref{gamma-id} reads, in particular, $\gamma-\frac{p-1}p=\frac{(1-\beta)(2-p)}p$. Moreover
\begin{equation}\label{eq:kappa}
\varkappa\ge-1\iff\varsigma\le\frac{N+1}{N-1}\iff p\ge p_* ,
\end{equation}
and $\varkappa>-1$ in our range.

\subsection{Ansatz and scaling}
For $u(x,t)=\tau^\beta F(x\tau^{-\gamma})$, $\tau=-t$, one has $\partial_tu=-\partial_\tau u=-\tau^{\beta-1}(\beta F-\gamma\,y\cdot DF)$ and $\Delta_pu=\tau^{(\beta-\gamma)(p-1)-\gamma}\,(\Delta_pF)(x\tau^{-\gamma})$. By \eqref{gamma-id}, $u$ solves $u_t=\Delta_pu$ if and only if $\Delta_pF=\gamma\,y\cdot DF-\beta F$. For a radial profile $F(y)=f(|y|)$, $f(0)=1$, and $g:=\phi_p(f')$, this is the system \eqref{ode}, whose integral form is the system \eqref{integral} of Section~\ref{section4}. This computation is valid for every $p>1$, but only where the profile is defined. Moreover, if $w$ solves $u_t=\Delta_pu$, so does $\lambda\,w(\nu x,t)$ whenever $\lambda^{p-1}\nu^p=\lambda$. With $\lambda=c$ and $\nu_c:=c^{\frac{2-p}p}$ this maps $\tau^\beta f(|x|\tau^{-\gamma})$ to $c\,\tau^\beta f(\nu_c|x|\tau^{-\gamma})$.

\begin{lem}[Local existence and uniqueness]\label{s:lem:loc}
For every $\beta\in\R$ there are $r_0>0$ and a \emph{unique} continuous solution $(f,g)$ of \eqref{integral} on $[0,r_0]$ with $|g(r)|\le M_\beta r$, $M_\beta:=\frac{(3/2)|\beta|+1}N$. It satisfies $f=1+O(r^{p'})$ and $g(r)=-\frac\beta Nr+O(r^{1+p'})$. In particular, for $\beta=0$ the solution is $f\equiv1$, $g\equiv0$.
\end{lem}
\begin{proof}
Let $\mathcal X:=\{g\in C([0,r_0]):\ |g(r)|\le M_\beta r\}$, with the sup norm, and let $\mathcal T(g)$ be the right-hand side of the first equation in \eqref{integral}, with $f=f_g:=1+\int_0^r\psi(g)$. For $g\in \mathcal X$ one has $|f_g(r)-1|\le\int_0^r(M_\beta s)^\varsigma ds\le\frac12$ if $r_0$ is small, hence
\[
|\mathcal T(g)(r)|\le r^{1-N}\!\int_0^r s^{N-1}\Big(|\gamma|M_\beta^\varsigma s^{1+\varsigma}+\tfrac32|\beta|\Big)ds\le\frac{\tfrac32|\beta|}Nr+\frac{|\gamma|M_\beta^\varsigma}{N+1+\varsigma}\,r^{1+\varsigma}\le M_\beta r
\]
for $r_0$ small, so $\mathcal T(\mathcal X)\subset \mathcal X$. Let $L:=\varsigma(M_\beta r_0)^{\varsigma-1}$, a Lipschitz constant of $\psi$ on $[-M_\beta r_0,M_\beta r_0]$, and $E:=\max|g_1-g_2|$. Then $|f_{g_1}-f_{g_2}|(s)\le LsE$ and
\[
|\mathcal T(g_1)-\mathcal T(g_2)|(r)\le r^{1-N}\!\int_0^r s^{N-1}\,(|\gamma|+|\beta|)\,L\,s\,E\,ds=\frac{(|\gamma|+|\beta|)L}{N+1}\,r^2E .
\]
Since $\varsigma>1$, $L\to0$ as $r_0\to0$, so $\mathcal T$ is a contraction on $\mathcal X$ for $r_0$ small, and Banach's theorem gives existence and uniqueness in $\mathcal X$. The expansions follow from \eqref{integral}: $f-1=\int\psi(g)=O(r^{\varsigma+1})$ and $\varsigma+1=p'$.

Uniqueness holds in fact among all continuous solutions of \eqref{integral}: any such solution has $|g(r)|\le Cr$ near $0$, because the integrand in \eqref{integral} is bounded, and the same contraction estimate, with $C$ in place of $M_\beta$, gives $E(r)\le\frac{(|\gamma|+|\beta|)L_C}{N+1}r^2E(r)$, hence $E\equiv0$ near $0$. For $r>0$ the system \eqref{ode} is an ODE with locally Lipschitz right-hand side ($\psi\in C^1$), so the solution has a unique maximal continuation to some $[0,r_{\max})$, $r_{\max}\in(0,\infty]$.
\end{proof}

\begin{prop}[Shape of the profile]\label{s:prop:shape}
Let $p_*<p<2$, and let $(f,g)$ be continued to its maximal interval $[0,r_{\max})$. Put $y_{\max}:=r_{\max}$.
\begin{enumerate}
\item[(a)] $\beta=0$: $f\equiv1$.
\item[(b)] $\beta<0$: $g>0$, $f\ge1$ is increasing, $y_{\max}=:y_b<\infty$ and $f(r)\to+\infty$ as $r\uparrow y_b$. In particular no profile with $\beta<0$ is global.
\item[(c)] $0<\beta<\beta_c$: $g<0$ on $(0,y_1]$, $f$ is strictly decreasing on $[0,y_1]$, where $y_1:=\sup\{r<r_{\max}:\ f>0\text{ on }[0,r)\}<r_{\max}$, $y_1<\infty$ and $f'(y_1)<0$.
\item[(d)] $\beta=1$: $f(r)=1-kr^{p'}$, $k=\frac{p-1}pN^{-\frac1{p-1}}$, $y_1=k^{-\frac1{p'}}$.
\end{enumerate}
\end{prop}
\begin{proof}
(a) is Lemma~\ref{s:lem:loc}. The sign of $g$ is proved by contradiction at a first zero $r_g$: there $(r^{N-1}g)'(r_g)=-\beta r_g^{N-1}f(r_g)$, which has the wrong sign as soon as $f(r_g)>0$ (for $\beta<0$ one has $f\ge1$). Let $H:=r^{N-1}|g|$.

(b) Here $\gamma>0$, $g>0$, and $H'=r^{N-1}(\gamma r g^\varsigma+|\beta|f)\ge\gamma r^\varkappa H^\varsigma$, because $r^{N-1}\,r\,g^\varsigma=r^\varkappa H^\varsigma$. Hence $(H^{1-\varsigma})'\le-(\varsigma-1)\gamma r^\varkappa$ and, for $r_*>0$,
\[
0<H(r)^{1-\varsigma}\le H(r_*)^{1-\varsigma}-(\varsigma-1)\gamma\int_{r_*}^rs^\varkappa ds .
\]
By \eqref{eq:kappa} $\varkappa>-1$, so the integral diverges and $y_b<\infty$. If $f$ were bounded on $[0,y_b)$ then $\int g^\varsigma=\int f'<\infty$, so $H'\le C(1+g^\varsigma)$ would be integrable and $H$, $g$ bounded, and the solution could be continued. So $f\to+\infty$.

(c) On $[0,y_1)$: $0<f\le1$, $G:=-g\ge0$, and $H'=r^{N-1}(\gamma rG^\varsigma+\beta f)>0$. Moreover $\int_0^rG^\varsigma=1-f(r)\le1$, so $H'\le C(1+G^\varsigma)$ and $H$ is bounded on bounded intervals; so is $G=r^{1-N}H$ away from $0$. Hence $y_1=r_{\max}<\infty$ is impossible. If $y_1=\infty$, then $f>0$ for all $r$ and $H'\ge\gamma r^\varkappa H^\varsigma$ yields, as in (b), a contradiction. So $y_1<r_{\max}$, $y_1<\infty$, and $H(y_1)>0$ gives $f'(y_1)<0$.

(d) $g=-r/N$ and $(r^{N-1}g)'=-r^{N-1}$ solve \eqref{ode} for $\beta=1$, $\gamma=\frac{p-1}p$, and uniqueness follows from Lemma~\ref{s:lem:loc}.
\end{proof}

\subsection{Positivity in $\Pa_h$}
For $\beta<\beta_c$ put $Y:=y_b$ if $\beta<0$, $Y:=\infty$ if $\beta=0$, $Y:=y_1$ if $\beta>0$, and, for $c>0$, $s\ge0$,
\[
u_{\beta,c,s}(x,t):=c(\tau+s)^\beta f\big(\nu_c|x|(\tau+s)^{-\gamma}\big),\quad\tau=-t,
\]
\[
y_c:=c^{\frac{p-2}p}Y,\quad\Phi_s(\tau):=R(\tau)(\tau+s)^{-\gamma}.
\]
By the ansatz and the scaling, $u_{\beta,c,s}$ is a weak solution wherever $\nu_c|x|(\tau+s)^{-\gamma}<y_{\max}$, and $u_{\beta,c,s}(0,-T)=c(T+s)^\beta$. Note that $Y\le y_{\max}$.

\begin{thm}\label{s:thm:positivity}
Let $p_*<p<2$, $\beta<\beta_c$, $c>0$, $s\ge0$, $h\ge0$. Then $u_{\beta,c,s}$ is defined and positive in $\Pa_h$ if and only if $\sup_{\tau>h}\Phi_s(\tau)\le y_c$. Consequently:
\begin{enumerate}
\item[(1)] $\beta=0$: $u\equiv c$.
\item[(2)] $\beta<1$, $\beta\neq0$, $s=0$: $u_{\beta,c,0}$ is defined and positive in $\Pa_h$ if and only if $h\ge h_Y(c):=h_Y(1)\,c^{\frac1{1-\beta}}$, where $h_Y(1)=\big(A^{-\frac{p-1}p}/Y\big)^{\frac p{(1-\beta)(2-p)}}$. For $s\ge0$ the condition $h+s\ge h_Y(c)$ is sufficient. The growth is $u(0,-T)=c(T+s)^\beta$.
\item[(3)] $\beta=1$: $u_{1,c,s}$ is positive in $\Pa_h$ if and only if $c\le a(A)$, where $a(A):=(A/k)^{\frac{p-1}{2-p}}$. Here $u_{1,c,0}=c\tau-c^{\frac1{p-1}}k|x|^{p'}$.
\item[(4)] $1<\beta<\beta_c$: $u_{\beta,c,s}$ is \emph{never} positive in $\Pa_h$.
\end{enumerate}
\end{thm}
\begin{proof}
\emph{Criterion.} A point $(x,t)$ lies in $\Pa_h$ iff $\tau>h$ and $|x|<R(\tau)$; then $z:=\nu_c|x|(\tau+s)^{-\gamma}<\nu_c\Phi_s(\tau)$, and $\nu_c\Phi_s\le Y$ for all $\tau>h$ if and only if $\sup\Phi_s\le y_c$. In that case $z<Y$: the profile is defined and, for $\beta>0$, $f(z)>0$ (for $\beta<0$, $f\ge1$). Conversely, if $\Phi_s(\tau)>y_c$ for some $\tau>h$, a point with $z=Y$ lies in $\Pa_h$: for $\beta>0$ there $u=0$, and for $\beta<0$ the solution is not defined there.

\emph{(2).} By \eqref{gamma-id}, $\gamma-\frac{p-1}p=\frac{(1-\beta)(2-p)}p>0$, so for $s=0$ the function $\Phi_0(\tau)=A^{-\frac{p-1}p}\tau^{-\frac{(1-\beta)(2-p)}p}$ is decreasing, with $\sup_{\tau>h}\Phi_0=A^{-\frac{p-1}p}h^{-\frac{(1-\beta)(2-p)}p}$, infinite if $h=0$. The condition $\sup\Phi_0\le c^{-\frac{2-p}p}Y$ is equivalent to $h\ge h_Y(1)c^{\frac1{1-\beta}}$. For $s\ge0$, $\Phi_s(\tau)\le A^{-\frac{p-1}p}(\tau+s)^{-\frac{(1-\beta)(2-p)}p}\le A^{-\frac{p-1}p}(h+s)^{-\frac{(1-\beta)(2-p)}p}$.

\emph{(3).} $\gamma=\frac{p-1}p$ and $\Phi_s(\tau)=A^{-\frac{p-1}p}\big(\frac\tau{\tau+s}\big)^{\frac{p-1}p}$, with supremum $A^{-\frac{p-1}p}$. The condition $A^{-\frac{p-1}p}\le c^{-\frac{2-p}p}y_1$, $y_1=k^{-\frac{p-1}p}$, is equivalent to $c^{2-p}\le(A/k)^{p-1}$, i.e. $c\le a(A)$.

\emph{(4).} $\gamma<\frac{p-1}p$, so $\Phi_s(\tau)\sim A^{-\frac{p-1}p}\tau^{\frac{(\beta-1)(2-p)}p}\to\infty$, whereas $y_c<\infty$ by Proposition~\ref{s:prop:shape}(c).
\end{proof}

\section{The singular supercritical case: the result in time}\label{section9}

We use the notation of Section~\ref{section5}; $R(\tau)$, $\delta_\mu$, $m_\mu$, $\Qt_\mu$ and $\zeta_\mu$ are defined as there, and are recalled below. The only change is the range of $\mu$, which is now $(0,\mu_1]$ with $\mu_1:=2^{-p'}\in(0,\frac14)$ instead of $(0,\frac14]$, because $\frac{p-1}p<\frac12$ when $p<2$ (it is exactly what makes $\delta_\mu R(\mu T)\le\frac12\delta_\mu R(T)$ in the proof of Proposition~\ref{s:prop:step2}). Let $\delta_\mu\in(0,1)$ be defined by $\delta_\mu^{p'}=\mu/2$ (it is a number depending on $\mu$ and $p$, not the constant $\delta$ of Proposition~\ref{exp pos sing}) and
\[
m_\mu(T):=\inf_{|x|<\delta_\mu R(T)}u(x,-T),\qquad\Qt_\mu(T):=\frac{m_\mu(T)}T,\qquad\zeta_\mu:=\frac{\ln2}{\ln(1/\mu)}\le\frac1{p'}<\frac12 .
\]

\begin{prop}[Persistence]\label{s:prop:persist}
There is $c_0\in(0,1)$ depending only on the data such that: if $Z=B_{2\rho}\times(t_0,t_0+S]$, $u\ge0$ is a continuous solution in a neighbourhood of $\overline Z$, $u(\cdot,t_0)\ge M$ in $B_{2\rho}$ and $\eps:=M^{p-2}S/\rho^p\le c_0$, then $u\ge M/2$ in $B_\rho\times(t_0,t_0+S]$.
\end{prop}
\begin{proof}[Differences with the proof of Proposition~\ref{prop:persist}]
For $p<2$ one has $\eps=S/(M^{2-p}\rho^p)$, and the DeGiorgi proof of Proposition~\ref{prop:persist} (Section~\ref{section5}) goes through with one change: the test function $w_n\chi_n^{p/2}$ must be replaced by $w_n\chi_n$, since $\chi_n^{\frac p2-1}$ is unbounded for $p<2$. Then $|D(w_n\chi_n)|^p\le2^{p-1}(\chi_n^p|Dw_n|^p+w_n^p|D\chi_n|^p)$, and $\displaystyle\int w_n^2\chi_n^2\le\int w_n^2\chi_n^p$ because $0\le\chi_n\le1$ and $p<2$. The exponent bookkeeping, which yields the factor $\eps^{\tilde\alpha}$, $\tilde\alpha=p/N$, is algebraic and does not depend on the sign of $p-2$.
\end{proof}
Clearly the proposition holds with any smaller $c_0$ as well. The threshold $q_\mu$ is given by the same expression as in \eqref{eq:qmu}, which for $p<2$ reads
\begin{equation}\label{s:eq:qmu2}
q_\mu:=\Big(\frac{2^{2p-1}}{c_0}\Big)^{\frac1{2-p}}\Big(\frac A\mu\Big)^{\frac{p-1}{2-p}} .
\end{equation}

\begin{prop}[Scale step]\label{s:prop:step2}
Let $\mu\in(0,\mu_1]$, $T>h/\mu$ and assume $u>0$ in $\Pa_h$. If $\Qt_\mu(T)\ge q_\mu$ then
\[
m_\mu(\mu T)\ge\tfrac12m_\mu(T)\qquad\text{and}\qquad\Qt_\mu(\mu T)\ge\frac1{2\mu}\,\Qt_\mu(T)\ \ \Big(\frac1{2\mu}\ge2^{p'-1}>2\Big).
\]
\end{prop}
\begin{proof}
Put $M:=m_\mu(T)$, $\rho:=\frac12\delta_\mu R(T)$, $t_0:=-T$, $S:=(1-\mu)T$, so that $B_{2\rho}=B_{\delta_\mu R(T)}$ and $t_0+S=-\mu T$.

\emph{The cylinder lies in $\Pa_h$.} For $\tau\in[\mu T,T]$: $\delta_\mu R(T)<\mu^{\frac{p-1}p}R(T)=R(\mu T)\le R(\tau)$ (because $\delta_\mu^{p'}=\mu/2<\mu$), and $\mu T>h$. So $\overline{B_{2\rho}}\times[-T,-\mu T]\subset\Pa_h$, $u>0$ there, and $u(\cdot,-T)\ge M$ in $B_{2\rho}$.

\emph{The condition $\eps\le c_0$.} We have $\rho^p=2^{-p}(\mu/2)^{p-1}(T/A)^{p-1}$, $M=\Qt_\mu(T)T$ and $S\le T$. Then
\[
\eps=\frac{M^{p-2}S}{\rho^p}\le\Qt_\mu(T)^{p-2}\,2^{2p-1}\mu^{1-p}A^{p-1}
\]
(the powers of $T$ cancel). Because $p-2<0$, the right-hand side is $\le c_0$ if and only if $\Qt_\mu(T)\ge q_\mu$.

Proposition~\ref{s:prop:persist} gives $u\ge M/2$ in $B_\rho\times(-T,-\mu T]$. At the time $-\mu T$ one has $\delta_\mu R(\mu T)=\mu^{\frac{p-1}p}\delta_\mu R(T)\le\frac12\delta_\mu R(T)=\rho$ iff $\mu\le2^{-p'}$, which holds. Hence $m_\mu(\mu T)\ge M/2$, and dividing by $\mu T$ gives the second claim.
\end{proof}

The hypothesis $\Qt\ge q_\mu$ is \emph{self-propagating}: since $\frac1{2\mu}>1$, the quotient does not decrease when one moves back in time.

\begin{cor}[Propagation]\label{s:cor:prop}
If $\tau_b>h$, $T\ge\tau_b$ and $\Qt_\mu(T)\ge q_\mu$, let $j\ge0$ with $T_j:=\mu^jT\in[\tau_b,\tau_b/\mu)$. Then $\Qt_\mu(T_i)\ge q_\mu$ for $i\le j$ and
\[
m_\mu(T_j)\ge2^{-j}m_\mu(T)=\Big(\frac{T_j}T\Big)^{\zeta_\mu}m_\mu(T).
\]
\end{cor}
\begin{proof}
Induction: if $\Qt_\mu(T_i)\ge q_\mu$ with $i<j$, then $T_i\ge T_{j-1}=T_j/\mu>h/\mu$, so Proposition~\ref{s:prop:step2} applies and gives $m_\mu(T_{i+1})\ge\frac12m_\mu(T_i)$ and $\Qt_\mu(T_{i+1})\ge\Qt_\mu(T_i)\ge q_\mu$. Finally $2^j=(\mu^{-j})^{\zeta_\mu}=(T/T_j)^{\zeta_\mu}$.
\end{proof}

\begin{thm}[A priori linear bound]\label{s:thm:time}
Let $u$ be a continuous weak solution of \eqref{p-Laplace}-\eqref{monotonicity} in an open set containing $\Pa_h$, with $u>\Theta$ in $\Pa_h$, and let $\mu\in(0,\mu_1]$; define $m_\mu$ with $u-\Theta$. Fix $\tau_b>h$ and put $C_b:=\max_{\tau\in[\tau_b,\tau_b/\mu]}(u(0,-\tau)-\Theta)$. Then for every $T\ge\tau_b$
\[
m_\mu(T)\le\max\Big\{q_\mu T,\ C_b\Big(\frac T{\tau_b}\Big)^{\zeta_\mu}\Big\}.
\]
In particular, with
\[
T_\mu:=\max\Big\{\tau_b,\ \big(C_b\tau_b^{-\zeta_\mu}/q_\mu\big)^{\frac1{1-\zeta_\mu}}\Big\},
\]
one has $\Qt_\mu(T)<q_\mu$ for all $T>T_\mu$. Therefore, there is \emph{no alternative}, and $\limsup_{T\to\infty}m_\mu(T)/T\le q_\mu$.
\end{thm}
\begin{proof}
Take $\Theta=0$ and $T\ge\tau_b$. If $\Qt_\mu(T)<q_\mu$ then $m_\mu(T)<q_\mu T$. Otherwise $\Qt_\mu(T)\ge q_\mu$, and by Corollary~\ref{s:cor:prop} and $m_\mu(T_j)\le u(0,-T_j)\le C_b$,
\[
q_\mu T\le m_\mu(T)\le\Big(\frac T{T_j}\Big)^{\zeta_\mu}m_\mu(T_j)\le C_b\Big(\frac T{\tau_b}\Big)^{\zeta_\mu}.
\]
This gives the inequality, and it is impossible when $T^{1-\zeta_\mu}>C_b\tau_b^{-\zeta_\mu}/q_\mu$, because $\zeta_\mu<1$. Hence $\Qt_\mu(T)<q_\mu$ for $T>T_\mu$.
\end{proof}

\begin{rem}[Consistency with the self-similar solutions]\label{s:rem:consistency}
Let $\beta=1$ and $c\le a(A)$. With $z_*^{p'}=\frac\mu{2A}c^{\frac{2-p}{p-1}}$ one finds $m_\mu(T)=cTf(z_*)$ and $\Qt_c=c\big(1-\frac\mu2\vartheta_c\big)\le c$, $\vartheta_c=\frac kAc^{\frac{2-p}{p-1}}\le1$: linear growth, with $\Qt_c\le a(A)$. Theorem~\ref{s:thm:time} requires $\Qt_c<q_\mu$ for every $c\le a(A)$ and every $\mu\le\mu_1$, which holds as soon as $a(A)\le q_{\mu_1}$, since $q_\mu$ decreases with $\mu$ and $\Qt_c<c\le a(A)$; by \eqref{s:eq:qmu2} and $k^{-(p-1)}=Np'^{\,p-1}$ this is $c_0\le2^{2p-1}\mu_1^{-(p-1)}k^{p-1}=2^{3p-1}k^{p-1}$. Since $c_0$ can always be decreased, this is only a consistency check on $c_0$. For $\beta<1$ one has $z\to0$ in the ball, so $m_\mu(T)\sim cT^\beta$ and $\Qt\to0$. For $1<\beta<\beta_c$ no positive solution exists (Theorem~\ref{s:thm:positivity}(4)).
\end{rem}

\section{Comparing the degenerate and the singular supercritical case}\label{section10}
In this section we summarize what the previous sections show about the two ranges. Statements for $p>2$ refer to Sections~\ref{section3}--\ref{section5}, those for $\frac{2N}{N+1}<p<2$ to Sections~\ref{section6}--\ref{section9}.

\subsection{Similarities}
\begin{rem}[Same tools, same geometry]
In both ranges the result in space uses only the expansion of positivity and the Comparison Principle (Theorems~\ref{P-L Thm in space} and~\ref{s:thm:space}), applied along a chain of cylinders whose vertices climb the paraboloid. The paraboloid $\Pa_h$ is the same set, its boundary is the zero set of the eternal solution of Remark~\ref{funzione paraboloidica}, and the natural growth in space is $|x|^{p'}$, $p'=\frac p{p-1}$, in both cases. In time, both rely on the persistence lemma (Propositions~\ref{prop:persist} and~\ref{s:prop:persist}), whose proofs differ only in the test function, and on the same scale step from $T$ to $\mu T$ on the ball $B_{\delta_\mu R(T)}$.
\end{rem}

\begin{rem}[Same profiles, same criterion]
The ansatz $u=\tau^\beta F(x\tau^{-\gamma})$ with $\gamma=\gamma_\beta$ leads, for every $p>1$, to the same ODE \eqref{ode}, to the same scaling $c\,\tau^\beta f(\nu_c|x|\tau^{-\gamma})$, and to the same positivity criterion: $u_{\beta,c,s}$ is positive in $\Pa_h$ if and only if $\sup_{\tau>h}\Phi_s(\tau)\le y_c$ (Theorems~\ref{thm:positivity} and~\ref{s:thm:positivity}). The case $\beta=1$ is explicit in both ranges, with $f(r)=1-kr^{p'}$, and the amplitude threshold $a(A)=(k/A)^{\frac{p-1}{p-2}}$ is the same expression in both cases.
\end{rem}

\begin{rem}[Same threshold in time]
The threshold $q_\mu$ of the scale step is given by the same expression \eqref{eq:qmu}, namely $q_\mu\simeq(\mu/A)^{\frac{p-1}{p-2}}$ up to a factor depending only on $c_0$ and $p$; the exponent $\zeta_\mu=\frac{\ln2}{\ln(1/\mu)}$ is the same, and in both cases the scale step compares $m_\mu(\mu T)$ with $m_\mu(T)$ through the factor $\frac12$.
\end{rem}

\subsection{Differences}

\begin{rem}[The common cause]
The intrinsic time scale of the expansion of positivity is $(\eta M)^{-(p-2)}\rho^p$ for $p>2$ (Proposition~\ref{exp pos}), which \emph{decreases} with $M$, and $\delta M^{2-p}\rho^p$ for $p<2$ (Proposition~\ref{exp pos sing}), which \emph{increases} with $M$. Likewise the smallness parameter $\eps=M^{p-2}S/\rho^p$ of the persistence lemma is small when $M$ is small if $p>2$, and when $M$ is large if $p<2$. All the reversals below come from this.
\end{rem}

\begin{rem}[Space: a lower bound on $\rho_0$ versus a lower bound on the slope]
For $p>2$ the construction in the proof of Theorem~\ref{P-L Thm in space} starts from a large radius, $\rho_0\ge\max\{|x_0|,(A\Lambda_\eta^{p'}u_0^{p-2}/\omega)^{\frac{p-1}{p(p-2)}}\}$ (see \eqref{eq:rho0}). In terms of the quotient $\Xi:=u_k/\rho_k^{p'}$, which is constant along the construction, the condition of Claim~1 is $A\Lambda_\eta^{p'}\le\omega\,\Xi^{2-p}$, that is $\Xi\le(\omega/A\Lambda_\eta^{p'})^{\frac1{p-2}}$: an upper bound on the slope, satisfied by taking $\rho_0$ large for any $u_0>0$. This is why the alternative of Theorem~\ref{P-L Thm in space} is ``either $u\le\Theta$ in $\{t<-h\}$, or $\limsup>0$''. For $p<2$ the same condition reads $\Xi\ge\Xi_*$ (Remark~\ref{s:rem:slope}): a lower bound on the slope. Hence $\rho_0$ has a maximum and not a minimum, it competes with $\rho_0\ge|x_0|$, and the alternative is the one of Theorem~\ref{s:thm:space}, in which the cone $\{u-\Theta<\Xi_*|x|^{p'}\}$ replaces the half-space. Moreover, the case distinction $M_k\ge2u_k$ or $M_k<2u_k$ in Claim~2 is not needed for $p<2$, since the radius $r_k\le\rho_k/16$ automatically.
\end{rem}

\begin{rem}[Explicit solutions $w_c$]
The computation of Proposition~\ref{s:prop:family} gives, for $p>2$, ``$w_c\le0$ in $\Pa_h$ if and only if $c\ge c_A$'', and for $p<2$ ``if and only if $c\le c_A$''. For $p>2$ the hypothesis $u\le\Theta$ in $\Pa_h$ forces a large slope $c$ at infinity within this family, in agreement with alternative (ii) of Theorem~\ref{P-L Thm in space}. For $p<2$ small slopes are admissible, which is why alternative (i) of Theorem~\ref{s:thm:space} cannot be improved to $u\le\Theta$ (Remark~\ref{s:rem:consequences}). The function of Remark~\ref{funzione paraboloidica} is $w_c$ with $\kappa_c=1$.
\end{rem}

\begin{rem}[Self-similar profiles]
\noindent
\begin{itemize}
\item For $p>2$ the function $\psi$ is not Lipschitz at $0$, and uniqueness for the integral system \eqref{integral} fails at $\beta=0$ (there are nonconstant solutions, see Theorem~\ref{thm:profile}); for $p<2$, $\psi$ is locally Lipschitz and Lemma~\ref{s:lem:loc} gives uniqueness for every $\beta$, by Banach's theorem.
\item For $p>2$ the profiles with $\beta<0$ are global, with $g>0$ and $y_1=\infty$ (Theorem~\ref{thm:profile}); for $p<2$ they blow up at a finite radius $y_b$ (Proposition~\ref{s:prop:shape}(b)), because $\psi$ is superlinear and $\varkappa>-1$. Hence for $p<2$ the corresponding solutions are defined only in $\{\nu_c|x|(\tau+s)^{-\gamma}<y_b\}$.
\item The sign in \eqref{gamma-id}, $\gamma-\frac{p-1}p=\frac{(\beta-1)(p-2)}p$, changes with $p-2$, and this exchanges the roles of the ranges of $\beta$ in the positivity theorems. For $p>2$ (Theorem~\ref{thm:positivity}): positive in $\Pa_h$ if and only if $\beta\le0$, or $\beta=1$ and $c\ge a(A)$, or $\beta>1$ and $h\ge h_0(c)$; never for $0<\beta<1$. For $p<2$ (Theorem~\ref{s:thm:positivity}): positive if and only if $\beta\le1$ (with $h\ge h_Y(c)$ for $\beta<1$, $\beta\ne0$, and $c\le a(A)$ for $\beta=1$); never for $1<\beta<\beta_c$. The inequality $c\gtrless a(A)$ is reversed.
\end{itemize}
\end{rem}

\begin{rem}[Time: a dichotomy versus an a priori bound]
The expression of $q_\mu$ is the same, $q_\mu=\big(c_02^{-(2p-1)}\mu^{p-1}A^{1-p}\big)^{\frac1{p-2}}$, but the exponent $\frac1{p-2}$ changes sign, and with it the inequality: for $p>2$ the persistence lemma applies when $\Qt_\mu(T)\le q_\mu$, for $p<2$ when $\Qt_\mu(T)\ge q_\mu$. Consequently:
\begin{itemize}
\item for $p>2$, $\Qt_\mu(T)\le q_\mu$ gives $m_\mu(T)\le2m_\mu(\mu T)$, and iterating one gets the dichotomy of Theorem~\ref{thm:dich}: either $\limsup\Qt_\mu\ge q_\mu$ (with $m_\mu\simeq T$ along a sequence) or $m_\mu(T)=O(T^{\zeta_\mu})$;
\item for $p<2$ the condition $\Qt_\mu\ge q_\mu$ propagates towards later times, and at a fixed time this is incompatible with $m_\mu\le u(0,\cdot)$, which is bounded on compact time intervals: Theorem~\ref{s:thm:time} is an a priori bound, with no alternative;
\item the restriction on $\mu$ is $\mu\le\frac14$ for $p>2$ (it uses $\frac{p-1}p>\frac12$), and $\mu\le2^{-p'}$ for $p<2$; accordingly $\zeta_\mu\le\frac12$ for $p>2$ and $\zeta_\mu\le\frac1{p'}$ for $p<2$.
\end{itemize}
In the consistency check with the profiles (Remarks~\ref{rem:consistency} and~\ref{s:rem:consistency}), for $p>2$ the profiles with $\beta=1$ and $\beta>1$ realise alternative (i) of Theorem~\ref{thm:dich}, with $\Qt\ge c(1-\frac\mu2)$ and $\Qt\to\infty$ respectively, and the profiles with $\beta\le0$ realise (ii); for $p<2$ all admissible profiles have bounded $\Qt$, in agreement with Theorem~\ref{s:thm:time}.
\end{rem}

\begin{rem}[From the infimum to the supremum]
For $p>2$ Theorem~\ref{thm:sup} turns the bound on the infimum into a linear bound on the supremum, through Lemma~\ref{lem:harnack}. That argument needs the forward Harnack cylinder, of height $\theta\rho^p$ with $\theta=(C_H/u_0)^{p-2}$, to fit in $\Pa_h$, which is the condition $u_0\ge c_1T$. For $p<2$ the same height is $(u_0/C_H)^{2-p}\rho^p$, and it fits if and only if $u_0\le c_1T$, that is, only for values which are already below the linear threshold. The step from the infimum to the supremum is therefore not transposed to the singular case.
\end{rem}


\section{Acknowledgments}
E. Z. wants to gratefully mention his current PhD supervisor at the University of Pavia, Ugo Gianazza. The initial interest in works as \cite{Gla}
and \cite{Shi} was aroused thanks to his mentorship. U. G. has then played a crucial role throughout the making of this work, with a number of significant discussions, decisive suggestions, and most useful corrections.
Also, E. Z. is a member of {\it GNAMPA} of {\it Istituto Nazionale per l’Alta Matematica (INdAM),} which partly supports his PhD program. 

\section{Disclaimer}
E. Z. states that Claude AI LLM has been used in the writing of this work for unsubstantial matters, such as grammar checks, language polishing, and typing in LaTeX.

\end{document}